\documentclass[12pt,english]{extarticle}

\usepackage[a4paper]{geometry}
\usepackage{fancyhdr}
\usepackage[dvipsnames]{xcolor}
\usepackage{latexsym,amssymb}
\usepackage{amsmath,amsthm,amsxtra,amsbsy,accents}
\usepackage{mathtools}
\usepackage{mathrsfs}
\usepackage{bm}
\usepackage{hyperref}

\DeclareMathAlphabet{\mathup}{OT1}{\familydefault}{m}{n}

\newcommand{\eps}{\varepsilon}
\DeclarePairedDelimiter{\abs}{\lvert}{\rvert}
\DeclarePairedDelimiter{\norm}{\lVert}{\rVert}
\DeclarePairedDelimiter{\bra}{(}{)}
\DeclarePairedDelimiter{\pra}{[}{]}
\DeclarePairedDelimiter{\set}{\{}{\}}
\DeclarePairedDelimiter{\skp}{\langle}{\rangle}
\makeatletter
\newcommand{\customlabel}[2]{%
   \protected@write \@auxout {}{\string \newlabel {#1}{{#2}{\thepage}{#2}{#1}{}} }%
   \hypertarget{#1}{}
}
\makeatother

\numberwithin{figure}{section}

\def\calD{{\mathcal D}} \def\calE{{\mathcal E}} 
  
  \def\calL{{\mathcal L}}
\def\calM{{\mathcal M}}  
  \def\calR{{\mathcal R}}
\def\calS{{\mathcal S}}  
  \def\calX{{\mathcal X}}
 
\def\rmd{{\mathrm d}}

 \def\rmB{{\mathrm B}} \def\rmC{{\mathrm C}}
\def\rmD{{\mathrm D}}  
 \def\rmH{{\mathrm H}} \def\rmI{{\mathrm I}}
  \def\rmL{{\mathrm L}}
\def\rmM{{\mathrm M}}

\def\rmV{{\mathrm V}} \def\rmW{{\mathrm W}} \def\rmX{{\mathrm X}}

  \def\scrL{{\mathscr  L}}
  
  \def\scrR{{\mathscr  R}}

\def\bbG{{\mathbb G}}  
  
 \def\bbN{{\mathbb N}} 
  \def\bbR{{\mathbb R}}
 \def\bbT{{\mathbb T}} 
  
 \def\bbZ{{\mathbb Z}}

\newtheorem{theorem}{Theorem}[section]
\newtheorem{lemma}[theorem]{Lemma}
\newtheorem{corollary}[theorem]{Corollary}
\newtheorem{definition}[theorem]{Definition}

\newtheorem{proposition}[theorem]{Proposition}
\newtheorem{remark}[theorem]{Remark}

\newtheorem{assumption}[theorem]{Assumption}

\makeatletter
\newcommand{\oset}[3][0ex]{%
  \mathrel{\mathop{#3}\limits^{
    \vbox to#1{\kern-2\ex@
    \hbox{$\scriptstyle#2$}\vss}}}}
\makeatother

\DeclareMathOperator*{\esssup}{ess\,sup}

\DeclareMathOperator{\shift}{S}

\DeclareMathOperator{\nablaDDiff}{{\overline{\bbG}_\swap}}
\DeclareMathOperator{\nablaDReact}{{\overline{\bbG}_\react}}
\DeclareMathOperator{\nablaDisc}{\overline{\bbG}}

\DeclareMathOperator{\nablaCont}{\bbG}
\DeclareMathOperator{\nablaCReact}{{\bbG_\react}}
\DeclareMathOperator{\nablaCDiff}{{\bbG_\swap}}

\DeclareMathOperator{\divDDiff}{{\overline{\bbG}^\ast_\swap}}
\DeclareMathOperator{\divDReact}{{\overline{\bbG}^\ast_\react}}

\DeclareMathOperator{\divCReact}{{\bbG^\ast_\react}}
\DeclareMathOperator{\divCDiff}{{\bbG^\ast_\swap}}

\newcommand{\N}{\mathbb{N}}
\newcommand{\R}{\mathbb{R}}

\newcommand{\cE}{\mathcal{E}}
\newcommand{\cR}{\mathcal{R}}

\newcommand{\Eu}{{E}}

\newcommand{\swap}{\mathrm{sw}}
\newcommand{\react}{\mathrm{re}}

\newcommand{\C}{\mathsf{C}}
\newcommand{\vertrule}[1][1.6ex]{\rule{.5pt}{#1}}
\newcommand{\CC}{\mathsf{C}\!\!\vertrule\;}

\newcommand{\AC}{\mathrm{AC}}

\newcommand{\dom}{\mathrm{dom}}
\newcommand{\tar}{\mathrm{tar}}

\newcommand{\discref}{w}
\newcommand{\mob}{m}

\newcommand{\babla}{\overline\nabla}

\newcommand{\dd}{\,\mathrm{d}}
\renewcommand{\div}{\mathop{\mathrm{div}}}

\definecolor{rot}{rgb}{1.000,0.000,0.000}

\newcommand{\one}{1\!\!1}

\newcommand{\CE}{\mathrm{CE}}
\newcommand{\lCE}{\overline{\mathrm{CE}}}

\usepackage{babel}

\newcommand{\ibar}{{\bar\imath}}
\newcommand{\idot}{{\vphantom{\imath}\smash{i}}}
\newcommand{\Ifull}{{I_\mathrm{full}}}

\newcommand{\Rone}{{R_\mathrm{full}}}
\newcommand{\Rtwo}{{R_\mathrm{lin}}}

\newcommand{\maxdensity}{M}

\newcommand{\CEn}{{C_\calE}}
\newcommand{\CDis}{{C_\calD}}

\newcommand{\CLone}{{C_{\rmL^1}}}
\newcommand{\CRate}{{C_\calR}}
\newcommand{\CSlope}{{C_\calS}}

\newcommand{\CBab}{{\overline{C}_{\babla}}}
\newcommand{\CBabsq}{{\widetilde{C}_{\babla}}}
\newcommand{\CProd}{{C_{\babla}}}

\usepackage{enumerate}

\newcommand{\grad}{\mathbb{G}}

\makeatletter

\renewcommand*\env@cases[1][1.2]{%
  \let\@ifnextchar\new@ifnextchar
  \left\lbrace
  \def\arraystretch{#1}%
  \array{@{\,}c@{\ }l@{}}%
}
\makeatother

\numberwithin{equation}{section}
\newcommand{\LB}{\lambda_{\mathrm B}}

\usepackage{authblk}
\author{Georg Heinze}

\affil{{
\small
Weierstrass Institute, Berlin, Germany, email:
\texttt{georg.heinze@wias-berlin.de}
}}

\title{Gradient flow formulation and discrete approximation of reaction-cross-diffusion systems with size exclusion}

\date{}

\begin{document}

\maketitle

\begin{abstract}
    \noindent 
    We characterize a class of possibly degenerate nonlinear reaction-cross-diffusion systems with size exclusion as (non-metric) gradient flows in continuity-equation format driven by entropy.
    We then establish a connection between these systems and discrete reaction networks by proving an EDP-type lower limit inequality for the associated gradient systems.
    In the non-degenerate case, we show a chain rule identity that allows us to identify minimizers of the corresponding energy-dissipation functional as weak solutions to the PDE system.
    In conjunction with the lower limit, this yields an alternative proof of existence of weak solutions and their approximability by solutions to the proposed discrete ODE systems.
    Finally, we generalize our results to models including drift terms by extending the established framework to tilted energy landscapes.

\end{abstract}

{
 \smallskip

 \smallskip

\noindent\emph{Acknowledgments:} The author is grateful to Alexander Mielke for stimulating discussions in particular about the setup of the discrete model. \smallskip

\noindent\emph{Funding:} This work was funded by the Deutsche Forschungsgemeinschaft (DFG, German Research Foundation) under Germany’s Excellence Strategy – The Berlin Mathematics Research Center MATH+ (EXC-2046/2, project ID: 390685689). \smallskip

\noindent\emph{AI Declaration:}
The author used Qwen3.8-Flash-Next and GPT‑5.5 Instant to search for literature, review elementary calculations/arguments after they had been carried out by the author, and proofread the manuscript.
The author takes full responsibility for the contents of the manuscript.
}

\section{Introduction}

In this work we are interested in a class of reaction-cross-diffusion systems with size exclusion,
\begin{equation}\label{eq:reaction-cross-diffusion_intro}
    \partial_t \rho_i 
    = \sum_{\ibar\in I}2K_{i\ibar} \div\bra[\big]{ 
    \rho_\ibar\nabla\rho_i - \rho_i\nabla\rho_\ibar}
	+ 
    \sum_{r\in R}\kappa_r \bra[\bigg]{
    \prod_{\ibar\in I}\rho_\ibar^{\alpha^r_\ibar}
    - \prod_{\ibar\in I}\rho_\ibar^{\beta^r_\ibar}}
    (\beta^r_i-\alpha^r_i),
\end{equation}
which we consider on the $d$-dimensional torus $\bbT^d$ for finitely many species $i\in I$ and finitely many mass-action-type reactions $r\in R$.
The cross-diffusion is characterized by symmetric non-negative coefficients $K_{i\ibar} = K_{\ibar i} \ge 0$, leading to non-uniform parabolicity whenever such a coefficient is zero.
Implicitly assuming detailed balance, the reactions are characterized by symmetrized reaction rates $\kappa_r>0$ and stoichiometric coefficients $\alpha^r_i,\beta^r_i \ge 0$ that conserve the total mass, i.e., satisfy $\sum_{i\in I} (\alpha^r_i -\beta^r_i) = 0$ for all $r\in R$.
Together with the symmetry $K_{i\ibar} = K_{\ibar i}$, this gives rise to the natural volume-filling constraint $\sum_{i\in I}\rho_i \equiv M$ for some $M>0$.

Our first goal is the derivation of a (non-metric) gradient system in continuity-equation format as introduced in~\cite{Peletier_Schlichting_2023}, which identifies~\eqref{eq:reaction-cross-diffusion_intro} as a gradient flow driven by the Boltzmann entropy and thereby encodes additional information about the underlying model beyond what is contained in the PDE alone.
Similar to~\cite{HMS25}, where reaction-diffusion systems without cross-diffusion were studied, the obtained gradient system contains non-convex primal and/or dual dissipation potentials, whenever reactions are nonlinear. 
This presents a challenging obstacle, which not only excludes Wasserstein-type metric structures (cf.~\cite{DNS09_CVPDE,Lisini2010,zinsl2015transport}), but also significantly complicates rigorous proofs of chain rule identities or convergence results.     
Depending on the coefficients $K_{i\ibar}$, the cross-diffusion adds to the lower order non-convexity of the reactions also a non-convexity in the highest order terms.

In the uniformly parabolic case, i.e., when $K_{i\ibar}>0$ for all pairs of species, we prove a chain rule identity, which in particular identifies weak solutions of~\eqref{eq:reaction-cross-diffusion_intro} with minimizers of the energy-dissipation functional of the proposed gradient system.
Note that this is more than what was obtained in~\cite{HMS25} for reaction-diffusion systems, where only a chain rule inequality was established, which is sufficient for approximation techniques, but generally does not identify the weak and variational solution concepts.

With the variational structure in place, our second goal is the recovery of~\eqref{eq:reaction-cross-diffusion_intro} from discrete mass-action-type reaction networks with detailed balance,
\begin{align}\label{eq:reaction_intro_general}
    \dot c_i 
    = \sum_{r\in \bar R}\bar\kappa_r \bra[\bigg]{
    \prod_{\ibar\in \bar I}c_\ibar^{\bar\alpha^r_\ibar}
    - \prod_{\ibar\in \bar I}c_\ibar^{\bar\beta^r_\ibar}}
    (\bar\beta^r_i-\bar\alpha^r_i),
\end{align}
for suitable reaction rates $\bar\kappa$ and suitable mass-conserving stoichiometric coefficients $\bar\alpha^r_i,\bar\beta^r_i \ge 0$.

More precisely, for each $N\in\bbN$ we denote $\bbZ^d_N = (\bbZ/N\bbZ)^d$ the set of discrete positions (with periodic boundary conditions) and $\Eu \coloneqq \set[\big]{e=(e_1,\ldots,e_d)^\top\in \set{0,1}^d, \sum_{l=1}^d e_l =1}$ the canonical Euclidean basis.
Then, we consider $\bar I = I{\times}\bbZ^d_N$ and $\bar R = (I^2{\times}\bbZ^d_N{\times}\Eu) \cup (R{\times}\bbZ^d_N)$ with the coefficients 
$\bar\kappa_{(i,\ibar,k,e)} = N^2 K_{i\ibar}$, 
$\bar\alpha^{(i,\ibar,k,e)}_{(\tilde\imath,\tilde k)} = \delta_{i,\tilde\imath}\delta_{k,\tilde k} + \delta_{\ibar,\tilde\imath}\delta_{k+e,\tilde k}$, 
$\bar\beta^{(i,\ibar,k,e)}_{(\tilde\imath,\tilde k)} = \delta_{i,\tilde\imath}\delta_{k+e,\tilde k} + \delta_{\ibar,\tilde\imath}\delta_{k,\tilde k}$ for the discrete cross-diffusion,
and $\bar\kappa_{(r,k)} = \kappa_r$, 
$\bar\alpha^{(r,k)}_{(i,k)} =\alpha^r_i$, 
$\bar\beta^{(r,k)}_{(i,k)} =\beta^r_i$ for the reactions, independent of the lattice site.
This turns~\eqref{eq:reaction_intro_general} into 
\begin{equation}\label{eq:reaction_intro}
\begin{aligned}
    \dot c_{i,k} 
    &= \sum_{\ibar\in I}\sum_{e\in \Eu} 2N^2 K_{i\ibar}\bra[\big]{
    c_{\ibar,k}\bra{c_{i,k+e} + c_{i,k-e}} - c_{i,k}\bra{c_{\ibar,k+e} + c_{\ibar,k-e}}
    }\\
    &+ \sum_{r\in R}\kappa_r \bra[\bigg]{
    \prod_{\ibar\in I}c_{\ibar,k}^{\alpha^r_\ibar}
    - \prod_{\ibar\in I}c_{\ibar,k}^{\beta^r_\ibar}}
    (\beta^r_i-\alpha^r_i).
\end{aligned}
\end{equation}

We interpret this system of ODEs as a gradient flow driven by the discrete Boltzmann entropy, thus also admitting a variational structure.
Allowing also for non-uniform parabolicity, i.e., $K_{i\ibar}= 0$ for some pairs of species (see assumptions~\eqref{as_general:pos_coeff_weak} and~\eqref{as_general:react}), we use purely variational arguments to show a lower limit inequality connecting the gradient systems inducing~\eqref{eq:reaction_intro} to the gradient system inducing~\eqref{eq:reaction-cross-diffusion_intro}. 
This lower limit is closely related to the notion of evolutionary $\Gamma$-convergence, which was first introduced in~\cite{SandierSerfaty,Serfaty2011}.
It connects the physical models underlying the respective equations and, when combined with the limit chain rule, which for solutions replaces a recovery sequence, in particular implies the convergence of solutions.

In the uniformly parabolic case, this yields an alternative proof of existence, though we emphasize that our main focus is not this result, but the identification and rigorous connection of the variational structures underlying both the discrete and continuum evolutions.

At the equation level the limit we consider is closely related to the convergence of finite volume schemes for cross-diffusion models.
These typically rely on a uniform $\rmH^{-1}$-control for time derivatives as established in~\cite{cances2019finite}.
In contrast, we derive compactness not only for solutions, but for all curves satisfying a discrete continuity equation and bounds on energy and dissipation.
Despite only having access to the weaker control provided by this dissipative structure, we obtain the necessary compactness by using different arguments in the near- and far-field regions, which we define in a non-trivial way.

Finally, we show that all the presented results, in particular also the approximation of gradient systems, remain true for models including drift potentials.
In the theory of evolutionary $\Gamma$-convergence for gradient
systems in the sense of the energy-dissipation principle (for short, EDP convergence), this is closely related to tilt-EDP convergence, a notion first introduced in~\cite{MiMoPe21EFED}.

\paragraph{Relation to literature.}

The system~\eqref{eq:reaction-cross-diffusion_intro} has several interesting special cases.
Indeed, choosing all $K_{i\ibar} = K/(2M) > 0$ for all $i,\ibar\in I$, we recover classical reaction-diffusion systems with shared diffusion constant and reactions of mass-action-type with detailed balance,
\begin{align*}
    \partial_t \rho_i 
    = K \Delta\rho_i 
    + \sum_{r\in R}\kappa_r \bra[\bigg]{
    \prod_{\ibar\in I}\rho_\ibar^{\alpha^r_\ibar}
    - \prod_{\ibar\in I}\rho_\ibar^{\beta^r_\ibar}}
    (\beta^r_i-\alpha^r_i).
\end{align*}
Using tools from ODE theory, it can be shown that such systems admit a unique classical solution, i.e., a solution satisfying a global $\rmL^\infty$-bound.
Classical solutions of reaction-diffusion systems have been extensively studied in the second half of the previous century (cf.~~\cite{Alik79LBSR, Rothe_1984, Morg89GESP, desvillettes2007global} and references therein), before it was discovered in~\cite{Pierre_Schmitt_1997} that these need not exist globally when species diffuse with different speeds, even for polynomial reactions.
This sparked an interest in weaker solution concepts, the most recent being renormalized solutions~\cite{Fischer, fischer2017weak}.
Furthermore, entropy methods have become popular in recent years, see e.g.\ in~\cite{desvillettes2006exponential, desvillettes2007global, desvillettes2014entropy, desvillettes2015entropic, DeFeTa17TERD}.
Finally, in~\cite{liero2013gradient, HMS25} reaction-drift-diffusion systems with mass-action-type kinetics and detailed balance were characterized as gradient flows that are driven by a relative Boltzmann entropy.

If instead of simplifying the cross-diffusion, we drop the reactions by setting $R = \emptyset$, we recover cross-diffusion systems with size exclusion,
\begin{align*}
    \partial_t \rho_i = \sum_{\ibar\in I} 2K_{i\ibar} \div\bra[\big]{ 
    \rho_\ibar\nabla\rho_i - \rho_i\nabla\rho_\ibar}.
\end{align*}
These models share with~\eqref{eq:reaction-cross-diffusion_intro} the non-uniform parabolicity when there are coefficients $K_{i\ibar}=0$.
They were popularized by~\cite{burger2010nonlinear} and have since been studied using both classical PDE as well as entropy methods, e.g.\ in~\cite{jungel2015boundedness, zamponi2017analysis, berendsen2017cross, Berendsen2019, chen2021cross, hopf2022multi}. 
In particular,~\cite{hopf2022multi} establishes a chain rule identity for the uniformly parabolic case.
We extend this result by also including reactions. 

System~\eqref{eq:reaction-cross-diffusion_intro} can also be viewed as a special case of the large class of reaction-diffusion systems studied in the uniformly elliptic case in~\cite{amann1990dynamic}, where local-in-time existence and uniqueness of classical solutions was established.
In particular, there are similarities to other reaction-cross-diffusion-type systems such as Maxwell-Stefan-type systems~\cite{maxwell1867dynamical,stefan1871}, which are used to describe multi-component fluids and, similar to the present model, also possess a volume-filling structure giving rise to
natural space-time-uniform $\rmL^\infty$-bounds.
Using entropy methods, such problems have recently been studied regarding global-in-time existence and long-time asymptotics, e.g.\ in \cite{jungel2013existence, daus2020exponential}.
Another large class of reaction-cross-diffusion systems is the Shigesada-Kawasaki-Teramoto (SKT) system originally suggested in~\cite{SKT1979} to study the spatial segregation of competing species in a heterogeneous environment. 
SKT-type systems usually do not possess a natural $\rmL^\infty$-control, but gain regularity from a linear diffusion term.
Here, we mention~\cite{chen2006analysis,desvillettes2014entropy,desvillettes2015entropic}, where entropy methods were used to show global-in-time existence of weak solutions for SKT-type systems with increasingly general cross-diffusion terms.

A main aspect of this work is the discrete-to-continuum limit, which we carry out on the level of variational structures, thus inherently preserving the volume-filling constraint and the variational description.
This is closely related to the approaches in~\cite{Schlichting2019} for the discrete-to-continuum limit of a growth model, in~\cite{esposito2025graph,esposito2023graph} for graph-to-local limits for single-species and multi-species nonlocal interaction equations, and in~\cite{lam2026convergence} for thermodynamic limits of stochastic particle systems.
Even more closely related are the finite-volume-type limits in~\cite{DisserLiero2015,HT23,heinze2024gradient,HeinzeKrautzPietschmann} for drift-diffusion-type equations with and without coupling to boundary evolutions in different geometries, in~\cite{HraivoronskaSchlichtingTse2024} for an aggregation-diffusion equation, in~\cite{mielke2025derivation} for the one-dimensional DLSS equation,
in~\cite{HMS25} for reaction-diffusion systems, and in~\cite{HMS26} for the porous-medium equation.
Furthermore, there are similarities with the convergence of entropy-structure- and volume-filling-constraint-preserving finite volume schemes, which have been studied in the context of cross-diffusion models~\cite{cances2020convergent}, reaction-cross-diffusion models for ion transport~\cite{cances2019finite,gaudeul2022entropy,cances2024convergence}, Maxwell-Stefan-type systems~\cite{cances2024finite}, and the SKT system~\cite{jüngel2021convergent}.

\paragraph{Structure of the paper.} 
In Section~\ref{s:abstract_strategy} we introduce the abstract notion of gradient system and gradient flow used in this work and relate them to other notions of gradient systems. 
Then, we formally sketch the strategy of approximating gradient systems.
In Section~\ref{s:GSRDS} we introduce the precise discrete and continuum models, discuss basic properties, and link them via constructing suitable embeddings for the discrete concentrations as well as their fluxes.
Section~\ref{s:MainResults} summarizes the assumptions and main results of this work.
The remaining sections are dedicated to proving the main results, first establishing compactness of the prelimit problem in Section~\ref{s:Compactness} and then a lower limit inequality in Section~\ref{s:Liminf}, thereby linking the discrete and continuum structures.
Section~\ref{s:CR} then establishes the equivalence of solution concepts and provides the remaining lower bound ensuring the convergence of discrete solutions to continuum solutions.
Finally, in Section~\ref{s:tilting} we extend the previous results to the tilted setting, where an additional drift potential is present in the evolution.

\section{Abstract strategy}\label{s:abstract_strategy}

To improve clarity, before challenging the reader with the notation of our concrete problem, we first present on a formal and abstract level the strategy of our work.

\subsection{Gradient systems with explicit abstract gradient mappings}
\label{subsec:GradSyst}

We begin by introducing a quintuple $(X,Y,\grad,\calE,\calR^\ast)$, called a \emph{gradient system}. This notion is a small modification of gradient systems in continuity-equation format introduced in \cite{Peletier_Schlichting_2023}.

The elements of the gradient system are two pairs of \emph{base spaces} 
$X = (X^\dom,X^\tar)$ and $Y = (Y^\dom,Y^\tar)$, where $X^\dom,Y^\dom$ are Borel 
subsets of a Euclidean space, and where $X^\tar,Y^\tar$ are Euclidean  spaces. 
Test functions over these spaces are linked by an abstract linear 
\emph{gradient map} $\grad : \rmC^\infty_c(X^\dom,X^\tar) \to \rmC^\infty_c(Y^\dom,Y^\tar)$,
which has a dual $\grad^\ast : (\rmC^\infty_c(Y^\dom,Y^\tar))^\ast \to \rmC^\infty_c((X^\dom,X^\tar))^\ast$, sometimes called 
\emph{process space to tangent space map}. Here,  for $Z\in\{X,Y\}$ the 
dual pairing is defined as usual by
\[
\forall\, \phi\in \rmC^\infty_c(Z^\dom,Z^\tar)\quad\forall\, \mu\in (\rmC^\infty_c(Z^\dom,Z^\tar))^\ast: \langle\phi,\mu\rangle_Z \coloneqq\langle\phi,\mu\rangle \coloneqq \int_{Z^\dom} \phi \cdot\dd\mu, 
\]
with $\cdot$ denoting the canonical inner product in the Euclidean space $Z^\tar$.

The fourth element of the quintuple is a lower semicontinuous (lsc) 
\emph{energy functional} $\calE:\calM_+(X^\dom,X^\tar)\to\R\cup\{\infty\}$, 
where $\calM_+(X^\dom,X^\tar)$ denotes the set of $X^\tar$-valued, 
component-wise non-negative  Radon measures. The final element is a 
\emph{dual dissipation potential} $\calR^\ast:\calM_+(X^\dom,X^\tar)\times 
C(Y^\dom,Y^\tar)\to [0,\infty]$, which, by definition, is lsc and non-negative 
 with $\calR^*(\rho,0)=0$,  and satisfies for all 
$\rho\in\calM_+(X^\dom,X^\tar)$ that $\xi\mapsto \calR^*(\rho,\xi)$ is convex.

Fixing an arbitrary time horizon $T>0$ and an initial datum 
$ \rho_0\in\calM_+(X^\dom,X^\tar)$, the gradient system $(X,Y,\grad,\calE,\calR^\ast)$
induces on $[0,T]$ an evolution equation, the \emph{gradient flow} equation 
\begin{subequations}
 \label{eq:abstract_GF}
\begin{align}
    \label{eq:abstract_GF_equation}
    &\partial_t \rho = \grad^\ast \partial_\xi\calR^\ast(\rho,-\grad \rmD\calE(\rho)), \\
    \label{eq:abstract_initial}
    &\rho(0) = \rho_0,
\end{align}    
\end{subequations}
where $\rmD\calE$ denotes the variational derivative of $\calE$ and $\partial_\xi\calR^\ast$ denotes the convex subdifferential of $\calR^\ast(\rho,\cdot)$. 

By specifying $\calR^\ast$ and $\calE$, the gradient system contains more information than the gradient-flow equation. Indeed, it is well-known that the same gradient-flow equation can be derived from different gradient systems, each corresponding to a different physical setting, see \cite{Mie24GF}. 

Before further discussing the link between gradient flow and gradient system, we comment on the relation of the presented notion of gradient flow with other notions.

\begin{remark}[Link to other notions of gradient flow]
\label{rem:link_to_other_GF-notions}
The presented notion is heavily influenced by the gradient systems in 
continuity-equation format introduced in \cite{Peletier_Schlichting_2023}, 
the only difference being the split of $X$ and $Y$ in a domain and a target 
space, which allows us to directly incorporate well-known objects like the 
classical gradient $\nabla: \rmC^1(\bbT^d)\to \rmC(\bbT^d;\R^d)$ into our framework.
Furthermore, we observe that by setting $\grad=\operatorname{id}$, $Y = X$, 
$\calX = \calM_+(X)$ and $\scrR^\ast(\rho,\varphi) \coloneqq 
\calR^\ast(\rho,\grad\varphi)$ for all $\rho,\varphi$, we recover the 
well-established notion of a gradient system $(\calX,\calE,\scrR^\ast)$ 
as introduced in \cite{Mie24GF}. Choosing $\calR^\ast$ as a quadratic functional, we can 
also recover metric gradient systems and metric gradient flows in the 
spirit of \cite{AGS}.
\end{remark}

To establish the link between the gradient system $(X,Y,\grad,\calE,\calR^\ast)$ and 
the gradient flow equation \eqref{eq:abstract_GF}, we split the latter into two parts:
First, we introduce the \emph{continuity equation}, which links a weak-$^\ast$ 
measurable curve $\rho:[0,T]\to \calM_+(X^\dom,X^\tar)$ 
with a  weak-$^\ast$ measurable \emph{curve of fluxes} $j:[0,T] \to \calM(Y^\dom,Y^\tar)$ 
by the relation (understood in the sense of distributions on $[0,T]\times X^\dom$)
\begin{align}\label{eq:abstract_CE}
    \partial_t \rho = \grad^\ast j.
\end{align}
The set of curves  $(\rho,j)$  satisfying \eqref{eq:abstract_CE} is denoted 
by $\CE$.

Secondly, given a pair $(\rho,j)\in\CE$, we recover \eqref{eq:abstract_GF} if the 
initial condition \eqref{eq:abstract_initial} holds and $j$ satisfies (in the 
sense of measures) the \emph{constitutive relation} 
\begin{align}\label{eq:abstract_KR}
    j = \partial_\xi\calR^\ast(\rho,-\grad \rmD\calE(\rho)).
\end{align}
One important  link between a gradient system and its induced gradient 
flow is called the \emph{energy-dissipation principle}. It is formally 
established as follows. We introduce the \emph{(primal) dissipation potential} 
$\calR:\calM_+(X^\dom,X^\tar) \times \calM(Y^\dom,Y^\tar)\to [0,\infty]$ as the 
convex dual of $\calR^\ast$ with respect to the second variable. Together 
$\CE$, $\calE$, $\calR$ and $\calR^\ast$  give rise to the \emph{dissipation functional}
\begin{align*}
    \calD(\rho,j)\coloneqq 
    \begin{cases}
        \int_0^T \calR(\rho,j) {+} \calR^\ast(\rho,-\grad\rmD\calE(\rho))\dd t & \quad\text{for } (\rho,j)\in\CE\\
        +\infty & \quad\text{for } (\rho,j)\notin\CE,
    \end{cases}
\end{align*}
and the \emph{energy-dissipation functional}
\begin{align*}
    \calL(\rho,j) \coloneqq \calE(\rho(T))-\calE(\rho(0))+\calD(\rho,j).
\end{align*}
 We say that  the gradient system satisfies the \emph{energy-dissipation
principle} if $\rho$ solving \eqref{eq:abstract_GF}  (in a suitable weak sense) 
is equivalent to $(\rho,j)$ solving $\CE$
and $\calL(\rho,j) = 0$. 

A crucial role in making this principle rigorous is played by the 
\emph{chain rule inequality} for the gradient system,  which means that 
$\calL(\rho,j)\ge 0$ holds for all $(\rho,j)\in\CE$. This name is  
motivated by the following formal calculation:
\begin{align*}
    \frac{\rmd}{\rmd t}\calE(\rho) &= \langle\rmD\calE(\rho),\partial_t\rho\rangle_X = -\langle -\grad\rmD\calE(\rho),j\rangle_Y \ge -\calR^\ast(\rho,-\grad\rmD\calE(\rho))-\calR(\rho,j),
\end{align*}
where the first equality is the classical chain rule, the second equality holds for $(\rho,j)\in\CE$, and the inequality follows from the duality of $\calR^\ast$ and $\calR$ (Young-Fenchel estimate). 
Integrating in time from 0 to $T$,  we obtain $\calL(\rho,j)\geq 0$. 

However, if $\calL(\rho,j)\le 0$ is imposed additionally, then we must have equality in the Young-Fenchel inequality for a.a.\ $t\in [0,T]$: 
\[
\langle -\grad\rmD\calE(\rho),j\rangle_Y  = \calR^\ast(\rho,-\grad\rmD\calE(\rho)) +\calR(\rho,j). 
\]
By the Fenchel equivalence, this implies that \eqref{eq:abstract_KR} holds a.e.\ on [0,T]. Plugging this into CE \eqref{eq:abstract_CE} shows that \eqref{eq:abstract_GF} holds. The opposite direction from \eqref{eq:abstract_GF} to $\calL(\rho,j)=0$ with $j$ from \eqref{eq:abstract_GF} is obvious.  

Of course, we will make these arguments rigorous for the reaction-diffusion systems under consideration.

\subsection{Convergence of gradient systems}
\label{subsec:GeneralConvergenceStrategy}

Having introduced abstract gradient systems and briefly discussed the 
energy-dissipation principle, we now want to discuss, on an abstract level, 
our strategy for obtaining the convergence of gradient flows. 

To this end, consider a family of approximating gradient systems $(X_N,Y_N, 
\grad_N,\calE_N,\calR^\ast_N)_{N\in\bbN}$ inducing $\C\calE_N$, $\calD_N$, and 
$\calL_N$ as before. As a first step, one proves that each 
$(X_N,Y_N,\grad_N,\calE_N,\calR^\ast_N)$ satisfies a chain rule inequality 
and an energy-dissipation principle. Next, one shows that for each $N\in\bbN$ 
and each initial datum $\rho_N^0$ with $\calE_N(\rho_N^0)<\infty$ there exists 
a solution pair $(\tilde\rho_N,\tilde\jmath_N)\in\C\calE_N$ with 
$\tilde\rho_N(0)=\rho_N^0$ and $\calL_N(\tilde\rho_N,\tilde\jmath_N)=0$.

Our aim is to connect the approximating gradient systems with a limit gradient 
system $(X,Y,\grad,\calE,\calR^\ast)$ inducing $\CE$, $\calD$, and $\calL$. 
For this, one has to show that $(X,Y,\grad,\calE,\calR^\ast)$ also satisfies 
chain rule inequality and energy-dissipation principle.

To establish the link, a candidate curve that might be a solution for the 
limit system needs to be obtained by applying a compactness argument to the family of prelimit solutions $(\tilde\rho_N,\tilde\jmath_N)_{N\in\bbN}$. However, the different gradient systems are defined 
over different base spaces, hence a unified space is needed in which compactness 
can be realized. To this end, one constructs an \emph{embedding operator} $\iota_N:\calM_+(X_N^\dom,X_N^\tar)\to \calM_+(X^\dom,X^\tar)$ and a \emph{discretization operator} $\iota_N^\ast:\rmC^\infty_c(X^\dom,X^\tar)\to \rmC^\infty_c(X_N^\dom,X_N^\tar)$
such that for all $\varphi\in \rmC^\infty_c(X^\dom,X^\tar)$ it holds
\begin{align*}
    \langle \iota_N\rho_N,\varphi\rangle_X = \langle \rho_N,\iota_N^\ast\varphi\rangle_{X_N}.
\end{align*}
For the fluxes one constructs a second embedding $\iota_{N,\textnormal{flux}}:\calM(Y_N^\dom,Y_N^\tar)\to \calM(Y^\dom,Y^\tar)$ such that for all $\varphi\in \rmC^\infty_c(X^\dom,X^\tar)$ it holds
\begin{align*}
    \langle \iota_{N,\textnormal{flux}}\,j_N,\grad\varphi\rangle_Y = \langle j_N,\grad_N\iota_N^\ast\varphi\rangle_{Y_N}.
\end{align*}
Since the continuity equation is understood in the sense of distributions, this implies that $(\rho_N,j_N)\in\C\calE_N$ if and only if $(\iota_N\rho_N,\iota_{N,\textnormal{flux}}\,j_N)\in\CE$. 

Now it is possible to prove that for each family $(\rho_N,j_N)_{N\in\bbN}$ with $(\rho_N,j_N)\in\C\calE_N$ and $\sup_{N\in\bbN}\sup_{t\in[0,T]}\calE_N(\rho_N(t))<\infty$ as well as $\sup_{N\in\bbN}\calD_N(\rho_N,j_N)<\infty$, there exists $(\rho,j)\in\CE$ with $\calE(\rho(0))<\infty$ and $\calD(\rho,j)<\infty$ such that (along a subsequence) $(\iota_N\rho_N,\iota_{N,\textnormal{flux}}\,j_N) \rightharpoonup^\ast (\rho,j)$. 
In particular, such a limit $(\tilde\rho,\tilde\jmath)$ exists for the family $(\tilde\rho_N,\tilde\jmath_N)_{N\in\bbN}$.

Next, one shows that for each family $(\rho_N,j_N)_{N\in\bbN}$ satisfying the a priori 
bounds as before and each limit $(\rho,j)$ of the embedded family, we have the liminf 
estimates 
\begin{align}\label{eq:abstract_lower_limit}
    \liminf_{N\to\infty}\calD_N(\rho_N,j_N) \ge \calD(\rho,j) 
    \qquad \text{and}\qquad \liminf_{N\to\infty}\calE_N(\rho_N(t)) 
    \ge \calE(\rho(t))\quad \forall\, t\in[0,T].
\end{align}
Notice that this inequality relates the dissipation functionals of the non-embedded
curves with the limiting dissipation functional of the limiting curve.
In particular, \eqref{eq:abstract_lower_limit} holds for the previously obtained 
family of solutions $(\tilde\rho_N,\tilde\jmath_N)_{N\in\bbN}$ 
and each of its limits $(\tilde\rho,\tilde\jmath)$.

To conclude that the limits are indeed solutions, we now assume the well-preparedness of initial data
\begin{align*}
    \iota_N\rho_N^0 \rightharpoonup^\ast \rho^0\qquad \text{and}\qquad\lim_{N\to\infty}\calE_N(\rho_N^0) = \calE(\rho^0) <\infty.
\end{align*}
The energy identity combined with \eqref{eq:abstract_lower_limit} and the limit 
chain rule inequality yield
\begin{align*}
    0 
    = \liminf_{N\to\infty}\calL_N(\tilde\rho_N,\tilde\jmath_N) 
    \ge \calL(\tilde\rho,\tilde\jmath) 
    \ge 0,
\end{align*}
from which the energy-dissipation principle of the limit gradient system implies that $\tilde\rho$ is solution starting at $\tilde\rho(0) = \rho^0$ and that $\tilde\jmath$ is given by the kinetic relation \eqref{eq:abstract_KR}.

\section{Reaction and reaction-cross-diffusion gradient systems}\label{s:GSRDS}
\label{s:Grad4RDS}

We want to describe the evolution of $i_*$ chemical species $\rmX_i$ with $i\in \{1,\dots, {i_*}\} \eqqcolon I $ undergoing cross-diffusion on the $d$-dimensional torus $\bbT^d$ and interacting according to $r_*$ chemical reactions
\[
\sum_{i\in I} \alpha^r_i \rmX_i \rightleftharpoons \sum_{i\in I} \beta^r_i \rmX_i, \quad r\in\{1,\dots,r_*\}\eqqcolon R.
\]
preserving the total mass, i.e., satisfying $\sum_{i\in I} \alpha^r_i =\sum_{i\in I} \beta^r_i$.

Throughout we assume that we have finitely many 
species and reactions,  i.e., $i_*,r_*\in \N $.  
In the following, we will also use the effective stoichiometric vectors  
$\gamma^r\eqqcolon\alpha^r {-} \beta^r \in \R^I$.  
Moreover, we fix reaction coefficients $\kappa_r>0$ (describing the reaction speed) 
for each reaction as well as cross-diffusion coefficients $K_{i\ibar}\ge0$ for each pair $i,\ibar \in I$.

The evolution can be written as the gradient flow of a gradient system, which we will introduce in this section. 
Before doing so, we first introduce the fully discrete gradient system that will later be used for the finite approximation. 
We then introduce the continuum gradient system, before linking the two to each other via an embedding.

\subsection{Reaction gradient system}

We present the gradient structure for the spatially discrete system with fixed $N\in\bbN$. 
Denoting $\bbZ^d_N = (\bbZ/N\bbZ)^d$ the set of discrete positions (with periodic boundary conditions), and $\Eu \coloneqq \set[\big]{e=(e_1,\ldots,e_d)^\top\in \set{0,1}^d, \sum_{l=1}^d e_l =1}$ the set of directions, we introduce the spaces 
\begin{align*}
    X_N&&&\coloneqq (X_N^\dom,X_N^\tar) &&\coloneqq (I\times\bbZ^d_N,\R),\\
    Y_{N,\swap}&&&\coloneqq (Y^\dom_{N,\swap},Y^\tar_{N,\swap})&&\coloneqq (I^2\times\bbZ^d_N\times\Eu,\R),\\
    Y_{N,\react}&&&\coloneqq (Y^\dom_{N,\react},Y^\tar_{N,\react})&&\coloneqq (R\times\bbZ^d_N,\R),\\
    Y_N&&&\coloneqq (Y^\dom_N,Y^\tar_N)&&\coloneqq (Y^\dom_{N,\swap}\times Y^\dom_{N,\react},Y^\tar_{N,\swap}\times Y^\tar_{N,\react}).
\end{align*}
We introduce the short notation $\rmC(X_N) \coloneqq \rmC(X_N^\dom;X_N^\tar)$ and analogously for all other spaces of functions/measures over these spaces. Furthermore, given a time interval $[0,T]$, we write $\rmC([0,T]\times X_N)\coloneqq \rmC([0,T]\times X_N^\dom;X_N^\tar)$ and use analogous notation for all other spaces of functions or measures over these base spaces.

Elements of the state space $\calM_+(X_N)$ are denoted by $c$ and will be called \emph{(discrete) chemical concentrations}. 

Abusing notation, we denote by $\langle\cdot,\cdot\rangle_N$ the dual products for vectors as well as components, e.g., for $(\xi,\zeta)\in \rmC(Y_N)$ and $(u,v)\in \calM(Y_N)$ we write
\begin{align*}
    \langle(\xi,\zeta),(u,v)\rangle_N &= \langle \xi,u\rangle_N + \langle \zeta,v \rangle_N 
    = \sum_{i\in I} \langle \xi_i,u_i \rangle_N + \sum_{r\in R}\langle \zeta_r, v_r \rangle_N \\
    &= \frac{1}{N^d}\sum_{k\in\bbZ^d_N}\bigg(\sum_{i,\ibar\in I}\sum_{e\in\Eu}  \xi_{i,\ibar,k,e} u_{i,\ibar,k,e} + \sum_{r\in R} \zeta_{r,k} v_{r,k}\bigg)
\end{align*}
and use analogous notations for other functions/measures defined over $X_N$ or $Y_N$.

Having fixed these base spaces, we now introduce the gradient operator, energy functional, and dual dissipation potential completing the gradient system.

\begin{definition}[Gradient operator]\label{def:disc_grad}
Given stoichiometric vectors $\alpha^r,\beta^r\in [0,\infty)^I$ satisfying $\sum_{i\in I} \alpha^r_i = \sum_{i\in I} \beta^r_i$, we recall the notation $\gamma^r \coloneqq\alpha^r{-}\beta^r$ and define the \emph{swap gradient operator} 
$\nablaDDiff$, the \emph{reaction gradient operator} $\nablaDReact$, and the \emph{(total) gradient operator} $\nablaDisc$ by
\begin{align*}
    \nablaDDiff&:\rmC^\infty_c(X_N) \to \rmC^\infty_c(Y_{N,\swap}),\quad 
    &&\nablaDDiff\varphi_{i,\ibar,k,e} \coloneqq (\babla\varphi_i)_{k,e} - (\babla\varphi_\ibar)_{k,e}\\
    & &&\phantom{\nablaDDiff\varphi_{i,\ibar,k,e} }\,\coloneqq (\varphi_{i,k+e} - \varphi_{i,k}) - (\varphi_{\ibar,k+e} - \varphi_{\ibar,k}),\\
    \nablaDReact&:\rmC^\infty_c(X_N) \to \rmC^\infty_c(Y_{N,\react}),\quad 
    &&\nablaDReact\varphi_{r,k} \coloneqq (\Gamma\varphi_k)_r
    \coloneqq \sum_{i\in I}\gamma^r_i\varphi_{i,k},\\
    \nablaDisc&:\rmC^\infty_c(X_N) \to \rmC^\infty_c(Y_N),\quad 
    &&\nablaDisc\varphi \coloneqq (\nablaDDiff\varphi, \nablaDReact\varphi),
\end{align*}
where we also introduced the \emph{discrete gradient} $\babla:\rmL^2(\bbZ^d_N)\to \rmL^2(\bbZ^d_N\times\Eu)$ and the \emph{stoichiometric matrix} (or \emph{Wegscheider matrix}) $\Gamma:\R^I\to\R^R$.

The dual operators corresponding to the above are given by
\begin{align*}
    \divDDiff&:(\rmC^\infty_c(Y_{N,\swap}))^* \to (\rmC^\infty_c(X_N))^*,
    &&\divDDiff F_{i,k} \coloneqq -\sum_{\ibar\in I} [(\overline\div F_{i,\ibar})_k-(\overline\div F_{\ibar,i})_k]\\
    & &&\phantom{\divDDiff F_{i,k}}\,\coloneqq \sum_{\ibar\in I} \sum_{e\in\Eu}[(F_{i,\ibar,k-e,e}{-}F_{i,\ibar,k,e}) - (F_{\ibar,i,k-e,e}{-}F_{\ibar,i,k,e})],\\
    \divDReact&:(\rmC^\infty_c(Y_{N,\react}))^* \to (\rmC^\infty_c(X_N))^*,
    &&\divDReact J_{i,k} \coloneqq (\Gamma^\ast J_k)_i = \sum_{r\in R}\gamma^r_iJ_{r,k},\\
    \nablaDisc^\ast&:(\rmC^\infty_c(Y_N))^* \to (\rmC^\infty_c(X_N))^*, 
    &&\nablaDisc^\ast(F,J) \coloneqq  \divDDiff F + \divDReact J,
\end{align*}
where we also introduced the \emph{discrete divergence} $\overline\div:\rmL^2(\bbZ^d_N\times\Eu)\to \rmL^2(\bbZ^d_N)$.
\end{definition}

    In accordance with the introduction, where we identified~\eqref{eq:reaction_intro} as a special case of~\eqref{eq:reaction_intro_general}, the swap gradient $\nablaDDiff$ can also be understood as a Wegscheider matrix with entries $\gamma^{(i,\ibar,k,e)}_{(\tilde\imath,\tilde k)} = \delta_{i,\tilde\imath}(\delta_{k+e,\tilde k}-\delta_{k,\tilde k}) - \delta_{\ibar,\tilde\imath}(\delta_{k+e,\tilde k}-\delta_{k,\tilde k})$ corresponding to binary reactions with stoichiometric coefficients that are either zero or one.

    Further, observe that since we assumed $\sum_{i\in I} \alpha^r_i = \sum_{i\in I} \beta^r_i$ for all $r\in R$, the kernel of $\nablaDisc$ contains functions that are the same for all species, i.e., $\varphi_i =\varphi$ for all $i\in I$. 

    Finally, we note that for $F_{i,\ibar} = -F_{\ibar,i} \in\calM(Y_{N,\swap})$, i.e., antisymmetric discrete diffusive fluxes, and $\varphi\in C(X_N)$, we have the identities
    \begin{align*}
        \sum_{i,\ibar\in I}\skp{\nablaDDiff\varphi_{i,\ibar},F_{i,\ibar}}_N
        = \sum_{i,\ibar\in I}\skp{\babla\varphi_i -\babla\varphi_\ibar,F_{i,\ibar}}_N
        = 2\sum_{i\in I}\skp[\bigg]{\babla\varphi_i,\sum_{\ibar\in I}F_{i,\ibar}}_N,
    \end{align*}
    highlighting another link between $\nablaDDiff$ and $\babla$.

Having defined the gradient operator, we now introduce the functionals characterizing our discrete gradient flow dynamics.
\begin{definition}[Free energy]\label{def:disc_energy}
We define the discrete \emph{free energy} as the (Boltzmann) entropy
\begin{equation}\label{eq:discrete_energy}
    \calE_N(c) \coloneqq \frac{1}{N^d}\sum_{i\in I}\sum_{k\in\bbZ_N^d}
    \LB \bra{c_{i,k}},
\end{equation}
where the Boltzmann function is defined by $\LB (r)\coloneqq r\log r-r+1$. 
\end{definition}
The variational derivative of $\calE_N(c)$ with respect to $c_i$ is given by
\begin{align*}
    \rmD_i \calE_N(c) = \log\bra{c_i}. 
\end{align*}

\begin{definition}[Dual dissipation potential]\label{disc_R*}
Given $K_{i\ibar} = K_{\ibar i} \ge 0$, $i,\ibar\in I$, 
we define the \emph{discrete dual dissipation potential} $\calR^\ast_N: \calM_+(X_N)\times \rmC(Y_N) 
\to[0, \infty)$ for $c \in\calM_+(X_N)$, 
$\xi \in \rmC(Y_{N,\swap})$, and $\zeta \in \rmC(Y_{N,\react})$ by
\begin{align*}
    \calR^\ast_{N, \swap} (c,\xi) &\coloneqq \frac{1}{N^d} \sum_{i,\ibar\in I}\sum_{k\in\bbZ_N^d} \sum_{e\in \Eu} N^2 \mob_{i,\ibar,k,e}(c)    \C^*(\xi_{i,\ibar,k,e}),\\
    \calR^\ast_{N, \react} (c,\zeta) &\coloneqq  \frac{1}{N^d} \sum_{r\in R}\sum_{k\in\bbZ_N^d}  \kappa_r \bra[\big]{c_k^{\alpha^r}c_k^{\beta^r}}^{1/2} \C^*(\zeta_{r,k}),\\
    \calR^\ast_N (c,(\xi,\zeta)) &\coloneqq \calR^\ast_{N, \swap} (c,\xi) + \calR^\ast_{N, \react} (c,\zeta),
\end{align*}
where we also introduced the smooth convex function $\C^*(r) \coloneqq 4\cosh(r/2)-4$, the mobility function $\mob_{i,\ibar,k,e}(c) \coloneqq K_{i\ibar} \sqrt{c_{i,k}c_{\ibar,k}c_{i,k+e}c_{\ibar,k+e}}$, and the shorthand notation $\rho^\alpha \coloneqq \prod_{i\in I} \rho_i^{\alpha_i}$.
\end{definition}

\begin{definition}[Gradient system and gradient flow]\label{def:disc_GS-GF}
The objects introduced above form the \emph{discrete gradient system}
$(X_N,Y_N,\nablaDisc,\calE_N,\calR^\ast_N)$. 
It is easy to check that, formally, the corresponding \emph{gradient flow equation}
\begin{align*}
	\dot c = \nablaDisc^\ast \partial_{(\xi,\zeta)}\calR^\ast_N(c,-\nablaDisc \rmD\calE_N(c))
\end{align*}
indeed is the discrete reaction-cross-diffusion system with size exclusion,
\eqref{eq:reaction_intro}.
\end{definition}

\begin{remark}\label{rem:discrete_grad-sys}
	We make the following additional observations:
	\begin{enumerate}[(i)]
		\item The discrete dissipation potential consists of two parts, $\calR^\ast_{N, \swap}$, which corresponds to two particles of different species exchanging their spatial positions, and $\calR^\ast_{N, \react}$, which corresponds to chemical reactions, i.e., particles changing their species at a fixed point in space. 
		\item The swap terms have the structure of binary reactions with the four involved stoichiometric coefficients being equal to $1$ and symmetrized reaction rate $N^2 K_{i\ibar}$.
		The scaling factor $N^2$ will later provide the transition to a quadratic dual dissipation in the limit $N\to \infty$. 
	\end{enumerate}
\end{remark}

Throughout, we will make use of  various  properties of the function $\C^*$
characterizing $\calR^\ast_N$ as well as its convex conjugate $\C$. We gather these 
properties in the following lemma:

\begin{lemma}\label{lem:PropertiesC}
    The convex function $\C^*:\R\to [0,\infty)$ defined by
    \begin{align*}
        \C^*(\zeta)=4 \cosh(\zeta/2)- 4
    \end{align*}
    has the convex conjugate 
    \begin{align*}
        \C(s) \coloneqq \sup_{\zeta \in \R} \{\zeta s - \C^*(\zeta)\} = 2s\mathop{\mathrm{arsinh}} ( s/2) - 2\sqrt{s^2{+}4}+4
    \end{align*}
    and satisfies for all $a,b>0$ the identities
    \begin{align} \label{eq:CStarandLog}
    \begin{aligned} 
         \sqrt{ab}\cdot\C^*(\log a- \log b) 
            &= 2\abs[\big]{\sqrt{a} - \sqrt{b}}^2
            \qquad\text{and}\qquad
            \sqrt{ab}\cdot(\C^*)'(\log a - \log b) = a - b.
    \end{aligned}
    \end{align}
\end{lemma}

\begin{proof}
These results are obtained by elementary calculations, see e.g.\ \cite[Lem.\,3.4]{HT23} for more details.
\end{proof}

 In addition to $\C$ we also need its so-called \emph{perspective function} 
$\CC:\bbR\times [0,\infty) \to [0,\infty]$, which is given by
\begin{equation*}
\CC(s|w)\coloneqq\sup_{\zeta\in\R}\{s\zeta- w \C^*(\zeta) \} = 
\begin{cases} w\C(s/w) &\text{for } w>0 , 
       \\ \chi_0(s) &\text{for }w=0. \end{cases}
\end{equation*}

For $c>0$ the dual dissipation potential induces a slope term by the relation 
$\calS_N(c)=\calR^\ast_N(c,-\nablaDisc\rmD \calE_N(c))$. This definition can then be extended to all 
$c\in\calM_+(X_N)$ (cf.\ \cite[Remark 3.7]{lam2024variational}) by exploiting the identity
\eqref{eq:CStarandLog} for $\C^*$ and $\log$. This yields the
\textit{relaxed slope}.

\begin{definition}[Relaxed slope]
\label{def:slope_disc}
The \textit{relaxed slope} $\calS_N:\calM_+(X_N)\to[0,\infty)$ is defined by
\begin{align*}
    \calS_{N,\swap}(c) &\coloneqq  \frac{1}{N^d} \sum_{k\in\bbZ_N^d} 
    \sum_{i,\ibar\in I}  \sum_{e\in \Eu} 2 N^2 K_{i\ibar} 
    \abs[\big]{\sqrt{c_{i,k+e}c_{\ibar,k}} - \sqrt{c_{i,k} c_{\ibar,k+e}}}^2,\\
    \calS_{N,\react}(c) &\coloneqq  \frac{1}{N^d} \sum_{k\in\bbZ_N^d} 
    \sum_{r\in R} 2\kappa_r 
    \big(c_k^{\alpha^r/2} 
    - c_k^{\beta^r/2}\big)^2,\\
    \calS_N(c) &\coloneqq \calS_{N,\swap}(c)+\calS_{N,\react}(c).
    \end{align*}
\end{definition}
We introduce the square roots $U_i {\coloneqq} \sqrt{c_i}$ and note that we may rewrite the swap part of the relaxed slope using the previously introduced discrete gradient $\babla$ as
\begin{align}
	\notag
    \calS_{N,\swap}(c) 
    &= \frac{1}{N^d} \sum_{k\in\bbZ_N^d} \sum_{i,\ibar\in I}  \sum_{e\in \Eu} 2N^2 K_{i\ibar} \abs[\big]{U_{i,k+e}U_{\ibar,k} - U_{i,k}U_{\ibar,k+e}}^2\\
    \notag
	&= \frac{1}{N^d} \sum_{k\in\bbZ_N^d} \sum_{i,\ibar\in I}  \sum_{e\in \Eu} 2N^2 K_{i\ibar} \abs[\big]{U_{\ibar,k} \babla U_{i,k,e} - U_{i,k} \babla U_{\ibar,k,e}}^2\\
	\label{eq:S_swap-a-mean}
    &= \frac{1}{N^d} \sum_{k\in\bbZ_N^d} \sum_{i,\ibar\in I} \sum_{e\in \Eu} 2N^2 K_{i\ibar} \abs[\big]{\tfrac12(U_{\ibar,k}+U_{\ibar,k+e}) \babla U_{i,k,e}  - \tfrac12(U_{i,k}+U_{i,k+e}) \babla U_{\ibar,k,e}}^2.
\end{align}
Next, we introduce the primal dissipation potential $\calR_N$.

\begin{definition}[Primal dissipation potential]
\label{def:R_disc}
We define the primal dissipation potential $\calR_N: \calM_+(X_N)\times \calM(Y_{N,\swap}) \times\calM(Y_{N,\react})\to[0,\infty)$ for a discrete concentration $c \in\calM_+(X_N)$, 
a discrete swap flux $F \in \calM(Y_N)$, and a reaction flux $J \in \calM(Y_N)$ by
\begin{align*}
	\calR_{N,\swap}(c,F) &\coloneqq \frac{1}{N^d}\sum_{k\in\bbZ_N^d} \sum_{i,\ibar\in I} \sum_{e\in \Eu}  \CC\left(F_{i,\ibar,k,e} \middle| N^2\mob_{i,\ibar,k,e}(c)\right),\\
    \calR_{N,\react}(c,J) &\coloneqq \frac{1}{N^d}\sum_{k\in\bbZ_N^d} \sum_{r\in R}\CC\left(J_{r,k} \middle| \kappa_r c_k^{(\alpha^r+ \beta^r)/2}\right),\\              
	\calR_N(c,F,J)&\coloneqq \calR_{N,\swap}(c,F) + \calR_{N,\react}(c,J),
\end{align*}
with the convention $0\cdot\infty = 0$.
\end{definition}
Note that $\calR_N$ is the convex conjugate of $\calR_N^\ast$ with respect to the second argument. 
Indeed, the dualities follow directly from the duality of 
$\C$ and $\C^\ast$.
Here, we see one advantage of using the explicit gradient operator, as it decouples the arguments of the different summands of $\calR^\ast_N$, and thus implies that its dual is again the sum of their individual duals, avoiding the calculation of an infimal convolution.

\begin{definition}[Continuity equation]\label{def:CE_disc}
    The operator $\nablaDisc$ gives rise to the continuity equation
    \begin{align*}
        \dot c = \divDDiff F + \divDReact J.
    \end{align*} 
    For $\maxdensity >0$, we denote $\lCE^\maxdensity_N$ the set of triples $(c,F,J)\in\AC([0,T];\calM_+(X_N))\times \rmL^1([0,T];\calM(Y_N))$ satisfying on $\bbZ^d_N$ the above equation and $\sum_{i\in I} c_i =\maxdensity$.
\end{definition}

\begin{remark}\label{rem:discrete_CE}
We make the following observations:
\begin{enumerate}[(i)]
    \item In the above definition the time independent constant $\maxdensity$ captures the additional invariance corresponding to the kernel of $\nablaDisc$ containing functions that coincide for all species.
    \item The domains of $X_N$ and $Y_N$ are finite sets, such that the topology for the measure spaces is irrelevant. 
    Only in the continuum case will it be important to use the weak-$^\ast$ topology of measures (cf.\ Definition \ref{def:CE_cont}).
    \item Following the proof of \cite[Lem.\,3.1]{Erbar2014}, we observe that the above definition is indeed well-posed for $F$ and $J$ satisfying $\rmL^1$-bounds in time.
\end{enumerate}
\end{remark}

The previously introduced objects together give rise to the following functionals.
\begin{definition}[Energy-dissipation functional]\label{def:EDF_disc}
    Fixing $M>0$, for $0\le s <t\le T$ we introduce the dissipation functional $\calD_N^{[s,t]}:\rmL^1([0,T];\calM_+(X_N)\times \calM(Y_{N,\swap})\times\calM(Y_{N,\react})) \to [0,\infty]$
    \begin{align}\label{eq:DiscreteD}
        \calD_N^{[s,t]}(c,F,J) \coloneqq 
        \begin{cases}
        \int_s^t [\calR_N(c(\tau),F(\tau),J(\tau)) + \calS_N(c(\tau))] \dd \tau & \quad\text{for } (c,F,J)\in\lCE^\maxdensity_N\\
        +\infty & \quad\text{otherwise}.
        \end{cases}
    \end{align}
    Furthermore, we introduce the energy-dissipation functional
    \begin{align*}
        \calL_N^{[s,t]}(c,F,J) \coloneqq \calE_N(c(t))- \calE_N(c(s)) + \calD_N^{[s,t]}(c,F,J).
    \end{align*}
    In the sequel we will often use the shorter notations $\calD_N(c,F,J) \coloneqq \calD_N^{[0,T]}(c,F,J)$ as well as $\calL_N(c,F,J) \coloneqq \calL_N^{[0,T]}(c,F,J)$.
\end{definition}

We conclude this section by specifying our notion of solution for  \eqref{eq:reaction_intro}.
For the subsequent analysis, it will be crucial to use the notion of \emph{energy-dissipation balance solutions} (in short \emph{EDB solutions}) that are based on the energy-dissipation functional $\calL_N$. 
Theorem \ref{thm:DiscEDB.ODE} will provide a rigorous connection between this notion and the ODE system \eqref{eq:reaction_intro} based on the corresponding chain rule as explained in Section \ref{subsec:GradSyst}.

\begin{definition}[Discrete EDB solutions]
\label{def:sol_disc}
We say $c\in \AC([0,T];\calM_+(X_N))$ is a \emph{discrete EDB solution} of 
\eqref{eq:reaction_intro} with initial datum $c_0\in \calM_+(X_N)$ with $\sum_{i\in I}(c_0)_{i,k} = M$ for all $k\in\bbZ^d_N$, if $c(0) = c_0$ and there exist fluxes $F,J$ such that $(c,F,J) \in \lCE^\maxdensity_N$ and for $0\leq s < t\leq T$ 
we have 
\begin{equation*}
\calL_N^{[s,t]}(c,F,J) =0.  
\end{equation*}
\end{definition}

The discrete system we study here can be rewritten as a special case of the system introduced in \cite[Section 3.1]{HMS25}. Hence, the existence of EDB solutions is an immediate consequence of \cite[Theorem~4.3]{HMS25} and standard ODE theory:
\begin{theorem}[Discrete EDB and \eqref{eq:reaction_intro}]
\label{thm:DiscEDB.ODE}
A function $ c \in \AC([0,T];\calM_+(X_N))$ is a solution to the discrete system \eqref{eq:reaction_intro} if and only if for $F$ and $J$ given by
\begin{align*}
    F_{i,\ibar,k,e} &= - N^2 K_{i\ibar} \bra{c_{\ibar,k}c_{i,k+e} - c_{i,k}c_{\ibar,k+e}}
    = - N^2 K_{i\ibar} \bra{c_{\ibar,k}\babla c_{i,k,e} - c_{i,k}\babla c_{\ibar,k,e}},\\ 
    J_{r,k} &= -\kappa_r \bra[\big]{c_k^{\alpha^r} - c_k^{\beta^r}},
\end{align*}
the triple $(c,F,J)$ is a discrete EDB solution in the sense of Definition \ref{def:sol_disc}.

In particular, such a solution exists.
\end{theorem}

\subsection{Reaction-cross-diffusion gradient system}

We present the gradient structure for the cross-diffusion system on the torus $\bbT^d$; as a short-hand we will use the name ``\emph{continuum system}''.   The base spaces we consider are  
\begin{align*}
    X&&&\coloneqq (X^\dom,X^\tar) &&\coloneqq (\bbT^d,\R^I),\\
    Y_\swap&&&\coloneqq (Y^\dom_\swap,Y^\tar_\swap)&&\coloneqq \big(\bbT^d,(\R^d)^{I^2}\big),\\
    Y_\react&&&\coloneqq (Y^\dom_\react,Y^\tar_\react)&&\coloneqq (\bbT^d,\R^R),\\
    Y&&&\coloneqq (Y^\dom,Y^\tar)&&\coloneqq (Y^\dom_\swap\times Y^\dom_\react,Y^\tar_\swap\times Y^\tar_\react).
\end{align*}
Furthermore, we recall the short notations $\rmC(X) \coloneqq \rmC(X^\dom;X^\tar)$ and $\rmC([0,T]\times X)\coloneqq \rmC([0,T]\times X^\dom;X^\tar)$ as well as analogous notations for all other spaces of functions/measures over these spaces.

Elements of the state space $\calM_+(X)$ are denoted by $\rho$ and will be called \emph{(continuum) chemical concentrations}. 
In our situation the measures will always have a density with respect to the Lebesgue measure on $\bbT^d$, which (slightly abusing notation) will also be denoted by $\rho$.

Similar to before, we abuse notation, denoting by $\langle\cdot,\cdot\rangle$ the dual products for vectors as well as components, e.g., for $(\xi,\zeta)\in \rmC(Y)$ and $(u,v)\in \calM(Y)$ we write
\begin{align*}
    \langle(\xi,\zeta),(u,v)\rangle &= \langle \xi,u\rangle + \langle \zeta,v\rangle 
    = \sum_{i,\ibar\in I} \langle \xi_{i,\ibar},u_{i,\ibar}\rangle + \sum_{r\in R}\langle \zeta_r,v_r\rangle \\
    &= \int_{\bbT^d}\bigg(\sum_{i,\ibar\in I}\sum_{l=1}^d \xi_{i,\ibar,l}(x)u_{i,\ibar,l}(x) + \sum_{r\in R}\zeta_r(x)v_r(x)\bigg)\dd x,
\end{align*}
and we use similar notations for other functions/measures defined over $X$ or $Y$.

\begin{definition}[Gradient operator]\label{cont_grad}
Given stoichiometric vectors $\alpha^r,\beta^r\in [0,\infty)^I$, 
$\gamma^r =\alpha^r-\beta^r$, $r\in R$, we denote by $\nabla$ the (classical) gradient 
with $\div$ the divergence, its negative dual. 
As before, we denote by $\Gamma$ the stoichiometric matrix. 
Then, the  continuum gradient operators are given by
\begin{align*}
    \nablaCDiff&:\rmC^\infty_c(X) \to \rmC^\infty_c(Y_\swap),\quad 
    &&\nablaCDiff\varphi_{i,\ibar}(x) \coloneqq \nabla\varphi_i(x)-\nabla\varphi_\ibar(x),\\
    \nablaCReact&:\rmC^\infty_c(X) \to \rmC^\infty_c(Y_\react),\quad 
    &&\nablaCReact\varphi_r(x) \coloneqq (\Gamma\varphi)_r = \sum_{i\in I}\gamma^r_i\varphi_i(x),\\
    \nablaCont&:\rmC^\infty_c(X) \to \rmC^\infty_c(Y),\quad 
    &&\nablaCont \varphi \coloneqq  (\nablaCDiff\varphi ,  
     \nablaCReact\varphi).
\end{align*}

Their dual operators are given by
\begin{align*}
    \divCDiff&:(\rmC^\infty_c(Y_\swap))^* \to (\rmC^\infty_c(X))^*,\quad 
    &&\divCDiff f_i(x) \coloneqq -\sum_{\ibar\in I} [\div f_{i,\ibar} - \div f_{\ibar,i}],\\
    \divCReact&:(\rmC^\infty_c(Y_\react))^* \to (\rmC^\infty_c(X))^*,\quad 
    &&\divCReact j_i(x) \coloneqq (\Gamma^\ast j)_i = \sum_{r\in R}\gamma^r_i j_r(x),\\
    \nablaCont^\ast&:(\rmC^\infty_c(Y))^* \to (\rmC^\infty_c(X))^*,\quad 
    &&\nablaCont^\ast(f,j) \coloneqq  \divCDiff f + \divCReact j.
\end{align*}
\end{definition}
As for the discrete case, we notice that $\ker\nablaCont$ contains functions coinciding for all species, $\varphi_i=\varphi$ for all $i\in I$.
Furthermore, we observe that for $f_{i,\ibar} = -f_{\ibar,i} \in \calM(Y_\swap)$ and $\varphi\in C^1(X)$, it holds
\begin{align*}
    \sum_{i,\ibar\in I}\skp{\nablaCDiff\varphi_{i,\ibar},f_{i,\ibar}} 
    = \sum_{i,\ibar\in I}\skp{\nabla\varphi_i -\nabla\varphi_\ibar,f_{i,\ibar}}
    = 2\sum_{i\in I}\skp[\bigg]{\nabla\varphi_i,\sum_{\ibar\in I}f_{i,\ibar}},
\end{align*}
highlighting another link between the swap gradient and the classical gradient. 
This identity is the reason for us keeping the explicit factor $2$ in~\eqref{eq:reaction-cross-diffusion_intro}.

As for the discrete system, we now introduce the functionals completing the continuum gradient system.
\begin{definition}[Free energy]\label{def:cont_energy}
The continuum \emph{free energy} functional is given as the entropy 
\begin{equation}\label{eq:torus_energy}
    \cE(\rho) =\sum_{i\in I}\int_{\bbT^d}\LB \bra{\rho_i}\dd x,
\end{equation}
where we recall the Boltzmann function $\LB (r)=r\log r-r+1$.
\end{definition}
The variational derivative of $\calE(\rho)$ with respect to $\rho_i$ is given by
\begin{align*}
    \rmD_i\cE(\rho) = \log\bra{\rho_i}. 
\end{align*}

\begin{definition}[Dual dissipation potential]\label{def:cont_R*}
The \emph{continuum dual dissipation potential} $\calR^*: \calM_+(X)\times \rmC(Y)\to[0,\infty)$ is then given for $\rho \in\calM_+(X)$, $\xi \in \rmC(Y_\swap)$, and $\zeta \in \rmC(Y_\react)$ by
\begin{align*}
    \calR^\ast_\swap (\rho,\xi) &\coloneqq \sum_{i,\ibar\in I}\frac{K_{i\ibar}}{2} \int_{\bbT^d} \rho_i\rho_\ibar \abs{\xi_{i,\ibar}}^2\dd x,\\
    \calR^\ast_\react (\rho,\zeta) &\coloneqq  \sum_{r\in R} \kappa_r \int_{\bbT^d}\bra[\big]{\rho^{\alpha^r}\rho^{\beta^r}}^{1/2} \C^*(\zeta_r)\dd x,\\
    \calR^*(\rho,(\xi,\zeta)) &\coloneqq \calR^\ast_\swap (\rho,\xi) + \calR^\ast_\react (\rho,\zeta).
\end{align*}
\end{definition}
As in the discrete setting, from now on we will write $\calR^*(\rho,\xi,\zeta)$ instead of $\calR^*\big(\rho,(\xi,\zeta)\big)$ and do the same for similar objects.

\begin{definition}[Gradient system and gradient flow]\label{def:cont_GS-GF}
The previously defined objects together form the gradient system $(X,Y,\nablaCont,\calE,\calR^\ast)$.

Similar to the discrete setting, a quick formal calculation shows (up to regularity) that the corresponding \emph{gradient flow equation}
\begin{align*}
	\partial_t \rho = \nablaCont^\ast \partial_{(\xi,\zeta)}\calR^*(\rho,-\nablaCont \rmD\calE(\rho))
\end{align*}
indeed is the continuum reaction-cross-diffusion system with exclusion \eqref{eq:reaction-cross-diffusion_intro}.

\end{definition}

For $\rho>0$ smooth the dual dissipation potential induces a slope term by the relation $\calS(\rho)=\calR^*(\rho,-\nablaCont\rmD \calE(\rho))$. This definition can then be extended to all $\rho\in\calM_+(X)$ (cf., e.g., \cite[(3.24)]{heinze2024gradient}) yielding the relaxed slope.

\begin{definition}[Relaxed slope]\label{def:slope_cont}
    The relaxed slope $\calS:\calM_+(X)\to[0,\infty)$ is defined by
    \begin{align*}
        \calS_\swap(\rho) &\coloneqq  \sum_{i,\ibar\in I} 2 K_{i\ibar} \int_{\bbT^d} \abs[\big]{\sqrt{\rho_\ibar}\nabla\sqrt{\rho_i}- \sqrt{\rho_i}\nabla\sqrt{\rho_\ibar}}^2\dd x,\\
        \calS_\react(\rho) &\coloneqq \sum_{r\in R} 2\kappa_r \int_{\bbT^d} \abs[\big]{\rho^{\alpha^r/2} - \rho^{\beta^r/2}}^2\dd x,\\
        \calS(\rho) &\coloneqq \calS_\swap(\rho)+\calS_\react(\rho).
    \end{align*}
\end{definition}

Next, we introduce the primal dissipation potential $\calR: \calM_+(X)\times \calM(Y)\to[0,\infty)$, which, as before, is given as the convex conjugate of $\calR^*$ with respect to the second argument.
\begin{definition}[Primal dissipation potential]\label{def:R_cont}
    The primal dissipation potential $\calR: \calM_+(X)\times \calM(Y_\swap)\times \calM(Y_\react)\to[0,\infty)$ is defined for a concentration measure $\rho \in\calM_+(X)$, a swap flux $f \in \calM(Y_\swap)$, and reaction flux $j \in \calM(Y_\react)$ by
    \begin{align*}
        \calR_\swap(\rho,f) &\coloneqq 
        \begin{cases}
            \sum_{i,\ibar\in I}\frac{1}{2K_{i\ibar} }\int_{\bbT^d}\frac{|f_{i,\ibar}|^2}{\rho_i\rho_\ibar}
            \dd x  &\text{for } f_{i,\ibar} = f_{i,\ibar}\dd x\ll \scrL^d,\\
            \infty & \text{otherwise},
        \end{cases}\\
        \calR_\react(\rho,j) &\coloneqq 
        \begin{cases}
            \sum_{r\in R} \int_{\bbT^d}\CC\left(j_r \middle|\kappa_r(\rho^{\alpha^r}\rho^{\beta^r})^{1/2}\right)\dd x 
             & \text{for }j_r = j_r\rmd x\ll \scrL^d,\\
            \infty & \text{otherwise},
        \end{cases}\\
        \calR(\rho,f,j)&\coloneqq 
        \calR_\swap(\rho,f) + \calR_\react(\rho,j),
    \end{align*}
    where we again made a slight abuse of notation using the same notation for measures and their Lebesgue densities, and where we used the convention $0\cdot\infty = 0$.
\end{definition}

As in the discrete case, the operator $\nablaCont$ gives rise to the continuity equation
\begin{align*}
    \partial_t \rho =  \nablaCont^\ast(f,j)= \divCDiff f + \nablaCReact^* j,
\end{align*}
which now generally has to be understood in the sense of distributions.
\begin{definition}[Continuity equation]\label{def:CE_cont}
    For $\maxdensity>0$, we denote $\CE^\maxdensity$ the set of triples 
    $(\rho,f,j)$ with $\rho\in\AC([0,T];\calM_+(X))\cap\rmL^\infty([0,T]\times X)$, $f\in\rmL^1([0,T];\calM(Y_\swap))$, and $j\in\rmL^1([0,T];\calM(Y_\react))$ 
    satisfying for all 
    $\varphi\in \rmC^1([0,T]\times X)$
    \begin{align}\label{eq:CE_distributional}
        \pra[\bigg]{\int_{\bbT^d}\sum_{i\in I}\varphi_i\dd \rho_i}_{t=0}^T        
        &= \int_0^T\bra[\bigg]{\int_{\bbT^d}\sum_{i\in I}\partial_t\varphi_i \dd \rho_i 
        + \int_{\bbT^d}\sum_{i,\ibar\in I} (\nabla\varphi_i - \nabla\varphi_\ibar) \cdot \dd f_{i,\ibar} 
        + \int_{\bbT^d}\sum_{r\in R}(\Gamma\varphi)_r \dd j_r}\dd t
    \end{align}
    and $\sum_{i\in I}\rho_i = \maxdensity$ $\scrL^d$-a.e.\ on $\bbT^d$. 
\end{definition}

\begin{remark}\label{rem:continuum_CE}
    We make the following remarks regarding the above definition:
    \begin{enumerate}[(i)]
        \item Analogously to the discrete setting, the sum of all the species being fixed to the time invariant constant $\maxdensity$ captures the additional invariance arising from the kernel of $\nablaCont$ containing test functions that coincide for all species.
        \item Recall that the absolute continuity in $\AC([0,T];\calM_+(X))$ has to be understood with respect to a metric generating the weak-$^\ast$ topology of $\calM_+(X)$.
        \item The notation $\rmL^1([0,T];\calM(Y_\swap))$ means weakly measurable functions s.t.\ $t\mapsto \|f(t)\|_{\calM(Y)}\in\rmL^1([0,T])$ (and the analogue holds for $j$).
        \item We mention for completeness that by \cite[Lem.~8.1.2]{AGS} and \cite[Lem.~3.1]{Erbar2014} for $f,j$ satisfying an $\rmL^1$-bound in time, which in particular holds when $\calD(\rho,f,j)<\infty$ (see below definition), the above distributional formulation can be strengthened to: For all $\varphi\in \rmC^1(X)$ and all times in $[0,T]$ it holds
    \begin{align*}
        \frac{\rmd}{\rmd t}\pra[\bigg]{\int_{\bbT^d}\sum_{i\in I}\varphi_i \dd \rho_i} 
        = \int_{\bbT^d}\sum_{i,\ibar\in I} (\nabla\varphi_i - \nabla\varphi_\ibar) \cdot \dd f_{i,\ibar} 
        + \int_{\bbT^d}\sum_{r\in R}(\Gamma\varphi)_r \dd j_r.
    \end{align*}
    \end{enumerate}
\end{remark}

The previously introduced objects give rise to the following functionals.
\begin{definition}[Energy-dissipation functional]\label{def:EDF_cont}
    Fixing $M>0$, for $0\le s < t \le T$ we introduce the dissipation functional $\calD:\rmL^1([0,T];\calM_+(X)\times \calM(Y))\to [0,\infty]$
    \begin{align}\label{eq:ContinuousD}
        \calD^{[s,t]}(\rho,f,j) \coloneqq 
        \begin{cases}
        \int_s^t [\calR(\rho(\tau),f(\tau),j(\tau)) + \calS(\rho(\tau))] \dd \tau & \quad\text{for } (\rho,f,j)\in\CE^\maxdensity\\
        +\infty & \quad\text{otherwise}.\end{cases}
    \end{align}
    Furthermore, we introduce the energy-dissipation functional
    \begin{align*}
        \calL^{[s,t]}(\rho,f,j) \coloneqq \calE(\rho(t))- \calE(\rho(s)) + \calD^{[s,t]}(\rho,f,j).
    \end{align*}
    We will often use the shorter notations $\calD(\rho,f,j) \coloneqq \calD^{[0,T]}(\rho,f,j)$ and $\calL(\rho,f,j) \coloneqq \calL^{[0,T]}(\rho,f,j)$.
\end{definition}

We conclude this section by specifying our notion of solutions for the continuum-space reaction-cross-diffusion gradient system \eqref{eq:reaction-cross-diffusion_intro}. 

\begin{definition}[Continuum EDI and EDB solutions for \eqref{eq:reaction-cross-diffusion_intro}]
\label{def:sol_cont}
We say $\rho\in \AC([0,T];\calM_+(X))$ is a \emph{continuum EDI solution} of \eqref{eq:reaction-cross-diffusion_intro} with initial datum $\rho_0\in \mathrm{dom}(\calE)$ with $\sum_{i\in I}(\rho_0)_i = M$ a.e., if $\rho(0) = \rho_0$ a.e., $\sup_{t\in [0,T]} \calE(\rho(t)) < \infty$, and there exist fluxes $f,j$ such that $(\rho,f,j) \in\CE^\maxdensity$,
$\calD(\rho,f,j)<\infty$, and for $0\leq s < t\leq T $ we have 
\begin{align*}
\calL^{[s,t]}(\rho,f,j)
 \le0. 
\end{align*}
We say it is a \emph{continuum EDB solution} of \eqref{eq:reaction-cross-diffusion_intro} if it is an EDI solution and in addition the triple $(\rho,f,j) \in\CE^\maxdensity$ satisfies
\begin{align*}
\calL^{[s,t]}(\rho,f,j) \ge 0. 
\end{align*}
\end{definition}

\subsection{Embedding}
\label{s:Embedding} 
As a crucial step for obtaining the convergence of gradient systems, we highlighted in Section~\ref{subsec:GeneralConvergenceStrategy} the construction of suitable embeddings connecting the prelimit spaces to the limit space.

\begin{definition}[Embedding and discretization operators]\label{def:const_embedding}
    We introduce for each $k\in\bbZ_N^d$ the cube 
    \begin{align*}
        Q^N_k \coloneqq \set{x\in\bbT^d: x_l \in [k_l/N,(k_l+1)/N), l=1,\ldots,d}.
    \end{align*}
    Next, we define the piecewise constant embedding operator $\iota_N:\calM(X_N)\to \calM(X)$ by setting
    \begin{align*}
        \iota_N (c_i)(x) &\coloneqq \rho_i(x) \coloneqq \sum_{k\in\bbZ_N^d} \one_{Q^N_k}(x) c_{i,k}.
    \end{align*} 
    Dual to it, we introduce the discretization operator $\iota^*_N:\rmC(X)\to \rmC(X_N)$ by setting
    \begin{align*}    
        \iota^*_N\xi_{i,k}\coloneqq(\iota^*_N\xi_i)_k&\coloneqq N^d\int_{Q^N_k} \xi_i\dd x.
    \end{align*}
    For the swap fluxes we introduce $\iota_{N,\swap}:\calM(Y_{N,\swap}) \to \calM(Y_\swap)$ defined by
    \begin{subequations}
    \label{eq:def.jotaFlux}
    \begin{equation}
    \label{eq:def.iotaFlux.diff}
         \iota_{N,\swap} F_{i,\ibar}(x)\coloneqq f_{i,\ibar}(x) = (f_{i,\ibar,e_1}(x),\ldots,f_{i,\ibar,e_d}(x))^\top
    \end{equation}
    where for $ e \in E$ we use the piecewise affine interpolation
    \begin{align}
    \label{eq:def.iotaFlux.diff.1D}
        f_{i,\ibar,e}(x)
        &\coloneqq \frac{1}{N}\sum_{k\in\bbZ_N^d}\bigg(\int_0^1\one_{Q^N_{k+\theta e}}(x)\dd \theta\bigg) F_{i,\ibar,k,e}\\
        &\,=\frac1N\sum_{k\in\bbZ^d_N} \one_{Q^N_k}(x)\bra[\big]{F_{i,\ibar,k-e,e}(1-(N x-k)\cdot e) + F_{i,\ibar,k,e}(N x-k)\cdot e }.\notag
    \end{align}
    Finally, for the reactive fluxes we again interpolate piecewise constantly, defining the respective embedding map $\iota_{N,\react}:\calM(Y_{N,\react})\to \calM(Y_\react)$ similar to above by
    \begin{equation}
    \label{eq:def.iotaFlux.react}
        \iota_{N,\react} J_r(x)\coloneqq j_r(x)\coloneqq\sum_{k\in\bbZ_N^d}\one_{Q^N_k}(x) J_{r,k}.
    \end{equation}        
    \end{subequations}
\end{definition}
    Note that using the above embedding operator we see that $\calE_N(c) = \mathcal{E}(\iota_N(c))$ and analogous identities hold for other integral functionals as well. 

Furthermore, by definition, for all $i,\ibar\in I$ and all $e\in\Eu$ we have the bounds
\begin{align*}
    \int_{\bbT^d}\abs{f_{i,\ibar,e}}(x)\dd x\le \frac{1}{N} \int_0^1\sum_{k\in\bbZ_N^d}\int_{\bbT^d}\one_{Q^N_{k+e\theta}}(x)
    \dd x\dd \theta \abs{F_{i,\ibar,k,e}} \leq 
    \frac{1}{N^d}\sum_{k\in\bbZ_N^d}\frac{\abs{F_{i,\ibar,k,e}}}{N},
\end{align*}
while for all $r\in R$ it holds
\begin{align*}
    \int_{\bbT^d}\abs{j_r}(x)\dd x\le \frac{1}{N^d}\sum_{k\in\bbZ_N^d}\abs{J_{r,k}}.
\end{align*}
        Denoting the discrete $\rmL^1$-norms on $\rmL^1(\bbZ^d_N)$ by 
$\|G\|_{\rmL^1_N}\coloneqq\frac{1}{N^d}\sum_{k\in\bbZ^d_N}|G_k|$, 
these estimates can be rewritten equivalently as
\begin{align}\label{eq:embedded-flux-bound}
\forall\, i,\ibar\in I,e\in \Eu:\quad \|f_{i,\ibar,e}\|_{\rmL^1}\leq \frac 1 N 
\|F_{i,\ibar,e}\|_{\rmL^1_N},\quad \forall\, r\in R:\quad \|j_r\|_1\leq \|J_r\|_{\rmL^1_N}.
\end{align}
To simplify notation, for all $N\in\bbN$ and $y\in\R^d$, we define the 
\emph{shift operator}
\begin{align}\label{eq:shift}
    \shift^N_y:\rmL^1(\bbT^d)\to \rmL^1(\bbT^d)
    \qquad\text{with}\qquad 
    \shift^N_y\phi(x)\coloneqq\phi\Big(\Big(x+\frac y N\Big)\!\!\!\!\!\mod 1\Big),
\end{align}
where the modulus is applied componentwise.

The following lemma uses the embeddings defined above to connect the discrete and continuum continuity equations  (recall $\lCE^\maxdensity_N$ from Definition \ref{def:CE_disc}).  

\begin{lemma}
\label{lem:CE.DiscCont}
    For each $c\in\calM_+(X_N)$, $\varphi\in \rmC(X)$, it holds \begin{align*}
        \langle\iota_N c,\varphi\rangle =\langle c,\iota_N^* \varphi\rangle_N.
    \end{align*}
    Furthermore, we have for $\varphi\in \rmC^1(X)$, $(F,J)\in\calM(Y_N)$ that
    \begin{align*}
        \langle\iota_{N,\swap} F,\nablaCDiff\varphi\rangle =\langle F,\nablaDDiff\iota_N^* \varphi\rangle_N \quad \text{and} \quad 
        \langle\iota_{N,\react} J,\nablaCReact\varphi\rangle =\langle J,\nablaDReact\iota_N^* \varphi\rangle_N.
    \end{align*}
    In particular, it holds $(c,F,J)\in\lCE^\maxdensity_N$ if and only if $(\iota_N c,\iota_{N,\swap}F,\iota_{N,\react}J) \in\CE^\maxdensity$.
\end{lemma}
\begin{proof}
For the first equality, we calculate for all $\xi\in \rmC(X)$
\begin{align*}
    \langle\xi_i,\iota_N c_i\rangle = \int_{\bbT^d} \xi_i \sum_{k\in\bbZ_N^d} c_{i,k}\one_{Q^N_k}\dd x = \frac{1}{N^d}\sum_{k\in\bbZ_N^d} \iota^*_N \xi_{i,k}  c_{i,k} = \langle\iota^*_N(\xi_i), c_i\rangle_N.
\end{align*}
Furthermore, we have for all $\varphi\in \rmC^1(X)$ that
\begin{align*}
    \langle F,\nablaDDiff\iota_N^\ast \varphi \rangle _N
    &= \frac{1}{N^d}\sum_{k\in\bbZ_N^d}\sum_{i,\ibar\in I}\sum_{e\in \Eu} F_{i,\ibar,k,e} 
    [(\iota_N^\ast \varphi)_{i,k+e}
    - (\iota_N^\ast \varphi)_{i,k}
    - (\iota_N^\ast \varphi)_{\ibar,k+e}
    + (\iota_N^\ast \varphi)_{\ibar,k}]\\
    &=\sum_{k\in\bbZ_N^d}\sum_{i,\ibar\in I}\sum_{e\in \Eu} F_{i,\ibar,k,e} \int_{Q^N_k} [\shift^N_e\varphi_i-\varphi_i-\shift^N_e\varphi_\ibar+\varphi_\ibar]\dd x.
\end{align*}
Rewriting for each $i\in I$ the integral by
\begin{align*}
    \shift^N_e\varphi_i(x)-\varphi_i(x) &= \int_0^1 \nabla\varphi_i\bigg(x+\frac{e}{N}\theta\bigg) \cdot \frac{e}{N}\dd \theta= \frac1N\int_0^1 \shift^N_{e\theta}\partial_{x_e}\varphi_i(x) \dd \theta,
\end{align*}
 and using Fubini for the integrable integrand, we thus obtain
\begin{align*}
    \langle F,\nablaDDiff\iota_N^\ast \varphi \rangle_N 
    &= \sum_{k\in\bbZ_N^d}\sum_{i,\ibar\in I}\sum_{e\in \Eu} \frac{F_{i,\ibar,k,e}}{N} \int_{Q^N_k}\int_0^1 \shift^N_{e\theta}\partial_{x_e}[\varphi_i - \varphi_\ibar](x) \dd \theta\dd x\\
    & =  \sum_{k\in\bbZ_N^d}\sum_{e\in \Eu} \frac{F_{i,\ibar,k,e}}{N} \int_{\bbT^d}\one_{Q^N_k}(x)\int_0^1 \shift^N_{e\theta}\partial_{x_e}[\varphi_i - \varphi_\ibar](x) \dd \theta\dd x\\   
    &  =  \sum_{k\in\bbZ_N^d}\sum_{e\in \Eu} \frac{F_{i,\ibar,k,e}}{N} \int_{\bbT^d}\int_0^1 \one_{Q^N_{k+e\theta}}(x)
    \partial_{x_e}[\varphi_i - \varphi_\ibar](x) \dd \theta\dd x\\
    &= \int_{\bbT^d}\sum_{e\in \Eu} f_{i,\ibar,e}(x) \partial_{x_e}[\varphi_i - \varphi_\ibar](x) \dd x =\langle f,\nablaCDiff \varphi\rangle .
\end{align*}
For the reactive flux, we simply observe that
\begin{align*}
\langle J, \nablaDReact \iota^*_N \varphi \rangle_N    
    &= \frac{1}{N^d} \sum_{k\in\bbZ_N^d}\sum_{i\in I} (\Gamma^TJ)_i N^d\int_{Q^N_k} \varphi_i\dd x 
    = \sum_{k\in\bbZ_N^d}\sum_{r\in R} J_{r,k} \int_{Q^N_k} \Gamma \varphi_r\dd x 
    =\langle j, \nablaCReact \varphi \rangle.
\end{align*}
In particular, for all $\varphi\in \rmC^1(X)$ we obtain 
$\langle \dot c, \iota_N^*\varphi\rangle_N = \langle\divDDiff F, \iota_N^*\varphi \rangle_N + \langle\divDReact J, \iota_N^*\varphi \rangle_N$ 
if and only if
$\langle \partial_t \iota_N c, \varphi\rangle = \langle \divCDiff f, \varphi \rangle + \langle\divCReact j,\varphi \rangle$.
 This concludes the proof. 
\end{proof}

\section{Main results}
\label{s:MainResults}

Before we state our main results, we fix the assumptions on our problem, which are used throughout the paper without always referring to them. 

\begin{assumption}[Assumptions on the coefficients]\label{ass:General}
    Throughout, we assume there exists a set of indices $\emptyset\ne\Ifull \subset I$ and positive constants $K_\ast,K^\ast\in(0,\infty)$, such that
    \begin{equation}
    \label{as_general:pos_coeff_weak}
    \tag{K}
    \begin{aligned}
        \forall\, i\in I,\ibar\in\Ifull:&\quad 0 < K_* \le K_{i\ibar} = K_{\ibar i} \le K^\ast <\infty,\\ 
        \forall\, i,\ibar\in I\setminus\Ifull:&\quad K_{i\ibar} = K_{\ibar i} = 0.
    \end{aligned}
    \end{equation}
    Only to prove the chain rule identity for the limit system, Theorem~\ref{thm:CRContinuous}, we will require the stronger assumption that $\Ifull = I$, i.e., the full connectivity
    \begin{align}
    \label{as_general:pos_coeff}
    \tag{K'}
        \forall\, i,\ibar\in I:\quad 0 < K_* \le K_{i\ibar} = K_{\ibar i} \le K^\ast <\infty.
    \end{align}
    Finally, only to prove the lower limit inequality in Theorem~\ref{thm:Comp_Liminf}, we assume in addition to the standing balance assumption $\sum_{i\in I}\alpha^r_i = \sum_{i\in I}\beta^r_i$ for all $r\in R$ that the set of reactions splits into two subsets whose stoichiometric coefficients satisfy one of two assumptions:
    \begin{align}\label{as_general:react}\tag{R}
        R = \Rone \cup \Rtwo
        \quad\text{with}\quad
        \begin{cases}
            \forall r\in \Rone\, \exists i,\ibar\in\Ifull: \alpha^r_i,\beta^r_\ibar>0,\\
            \forall r\in \Rtwo\, \exists i,\ibar\in I: \alpha^r_{\tilde\imath} = \delta_{i,{\tilde\imath}}, \beta^r_{\tilde\imath} = \delta_{\ibar,{\tilde\imath}},
        \end{cases}
    \end{align}
    where by $\delta$ we denoted the Kronecker delta.
\end{assumption}

\begin{remark}
    Let us briefly discuss these assumptions:
    \begin{enumerate}[(i)]
    \item Heuristically, the assumption \eqref{as_general:pos_coeff} corresponds to all species being able to swap with each other, while \eqref{as_general:pos_coeff_weak} with $\Ifull=\set{0}$ is the natural choice in situations where the species $i=0$ embodies vacancies, since the jump of a particle into an unoccupied space can also be interpreted as the particle changing places with a vacancy.
    
    \item The authors of \cite{hopf2022multi} additionally consider cross-diffusion systems with size exclusion with ``partially interacting species", i.e., the situation where a species can swap places with some, but not all other species.
    However, in this situation they only show a weakened version of an energy-dissipation inequality, where they set both the slope and the dissipation to zero on sets, where all fully connected species vanish (cf.\ \cite[Eq. (4.4)]{hopf2022multi}) and an improvement of that result to a dual structure is not obvious.
    Therefore, we exclude this case in our considerations, only noting that our discrete prelimit model is also well-posed in that partially degenerate case.

    \item The chain rule identity for the degenerate limit system remains an open question.
    Indeed, to the best of the author's knowledge, even for systems without reactions the chain rule is not known in the degenerate setting.

    \item Regarding the assumption~\eqref{as_general:react}, first note that it is automatically implied by~\eqref{as_general:pos_coeff}.
    Indeed, in this case $R=\Rone$, which will allow us to obtain strong $\rmL^1$-convergence of $(\rho^N)^\lambda \to \rho^\lambda$ for all multiindices $\lambda\in\bbR^I_{\ge 0}$.
    In the more general situation, this strong convergence may fail on the set $\Rtwo$, thus leaving us with only weak convergence. 
    In this case, we rely on convexity arguments to deduce lower limit inequalities, which are valid simultaneously for both the rate term $\calR_\react$ and the slope term $\calS_\react$, only if the reactions are linear.
    \end{enumerate}
\end{remark}
Existence of discrete EDB solutions is known from Theorem~\ref{thm:DiscEDB.ODE}. 
Our goal is to use such solutions to approximate the continuum system and ensure that their limits are EDB solutions of the limit problem. 
However, a standard Aubin-Lions-type argument to obtain compactness is not applicable under the more general assumption \eqref{as_general:pos_coeff_weak}.
Instead, we adapt Kolmogorov-Riesz-type arguments to the inhomogeneous cosh-type gradient structure. We stress that for this we do not rely on the $H^1$-$H^{-1}$-type duality used in \cite[Lemma 10]{cances2019finite}, but only make use of the weaker control provided by the dissipation functional.
\begin{theorem}[Compactness and lower limit inequalities]
\label{thm:Comp_Liminf}
Assume~\eqref{as_general:pos_coeff_weak},~\eqref{as_general:react} and let $T>0$.
For $N\in\bbN$ consider
$(c^N,F^N,J^N)\in\lCE^\maxdensity_N$ satisfying the bounds
$\sup_{N\in\N}\esssup_{t\in[0,T]}\calE_N(c^N(t)) 
<\infty$ and $\sup_{N\in\N}\calD_N(c^N,F^N,J^N) <\infty$.

Then, there exists $(\rho,f,j)\in\CE^\maxdensity$ with $f\in\rmL^1([0,T]\times Y_\swap)$,
$j\in\rmL^1([0,T]\times Y_\react)$, and $\rho\in\rmL^\infty([0,T]\times X)$ satisfying $\sqrt{\rho_\ibar} \in\rmL^2([0,T];\rmH^1(\bbT^d))$ and $\sqrt{\rho_i} \in\rmL^2([0,T];\rmH^1(\bbT^d,\rho_\ibar))$ for all $i\in I$, $\ibar\in\Ifull$.
For this triple (up to a subsequence) we have 
$\iota_{N,\swap}F^N\rightharpoonup f$ weakly in 
$\rmL^1([0,T]\times Y_\swap)$, and $\iota_{N,\react}J^N 
\rightharpoonup j$ weakly in $\rmL^1([0,T]\times Y_\react)$, and $\iota_Nc^N_\ibar \to \rho_\ibar$ as well as $\iota_N(c^N_\ibar c^N_i)\to \rho_\ibar\rho_i$  strongly in $\rmL^p([0,T]\times \bbT^d)$ for all $i\in I$, $\ibar\in\Ifull$, $p\in[1,\infty)$.

Moreover, we have the lower limit inequalities
\[
    \liminf_{N\to\infty} \calD_N(c^N,F^N,J^N)\geq \calD(\rho,f,j) \quad 
    \text{and }\quad \liminf_{N\to\infty} \calE_N(c^N(t))\geq \calE(\rho(t))
     \text{ for all } t\in[0,T],
\]
for the functionals defined in \eqref{eq:discrete_energy}, \eqref{eq:DiscreteD},
\eqref{eq:torus_energy}, and \eqref{eq:ContinuousD}, respectively.

In particular, for well-prepared initial data, i.e.,
$\iota_N c^N(0)\to \rho(0)$ with 
$\calE_N(c^N(0))\to \calE(\rho(0))$, it holds
\begin{align*}
    \liminf_{N\to\infty} \calL_N(c^N,F^N,J^N)\geq \calL(\rho,f,j).
\end{align*} 
\end{theorem}

To conclude that the limit $(\rho,f,j)$ solves the limit gradient-flow equation, we additionally need a chain rule identity for the continuum 
system. 

Under the stronger assumption \eqref{as_general:pos_coeff}, we achieve this despite neither the rate nor the slope integral functionals possessing (jointly) convex integrands by adapting the technique from \cite{HMS25} and relying on the exclusion dynamics' natural $\rmL^\infty$-bounds.
\begin{theorem}[Chain rule identity]
\label{thm:CRContinuous}
Assume~\eqref{as_general:pos_coeff} and let $T>0$.
Consider a curve $(\rho,f,j)\in\CE^\maxdensity$ with 
$\sup_{t\in[0,T]}\calE(\rho(t)) <\infty$, and $\calD(\rho,f,j) <\infty$. 
Then, for every $0\le s < t\le T$ it holds
\begin{align*}
    \calE(\rho(t)) - \calE(\rho(s))
    &= \int_s^t\int_{\bbT^d}\bra[\Bigg]{\sum_{i,\ibar\in I} \bra[\big]{\nabla\rmD_i\calE(\rho) - \nabla\rmD_\ibar\calE(\rho)} \cdot f_{i,\ibar} 
    + \sum_{r\in R}(\Gamma\rmD\calE(\rho))_r j_r} \dd x \dd \tau.
\end{align*}
In particular, it holds  $\calL(\rho,f,j) \ge 0$ and $\rho$ is 
a continuum EDB solution of \eqref{eq:reaction-cross-diffusion_intro} in the sense of Definition~\ref{def:sol_cont} for fluxes respectively given by the constitutive relations
\begin{equation} 
\label{eq:const_rel_cont}
    f_{i,\ibar} = -K_{i\ibar} \bra{\rho_\ibar\nabla\rho_i -\rho_i\nabla\rho_\ibar}
    \quad \text{and} \quad
    j_r = -\kappa_r \bra[\big]{\rho^{\alpha^r} - \rho^{\beta^r}}
    \quad \text{a.e.\ on $[0,T]\times\bbT^d$ $\forall$ $i,\ibar\in I$, $r\in R$ }
\end{equation}
if and only if $\rho$ is a weak solution of the PDE~\eqref{eq:reaction-cross-diffusion_intro}.

\end{theorem}

Together with Theorem~\ref{thm:DiscEDB.ODE}, these theorems imply that solutions of the discrete problems \eqref{eq:reaction_intro} on $\bbZ^d_N$ starting from well-prepared initial data converge (after choosing a suitable subsequence) to solutions of the cross-diffusion system with exclusion \eqref{eq:reaction-cross-diffusion_intro}. 
This is summarized in the following corollary.

\begin{corollary}[Convergence of solutions]\label{cor:conv_sol}
    Assume~\eqref{as_general:pos_coeff} and let $T>0$.
    Let $\rho_0\in \rmL^1(X)$ s.t.\ $\rho_0\ge 0$ and $\sum_{i\in I}\rho_{0,i}=\maxdensity$ a.e.\ on $\bbT^d$, and $\calE(\rho_0)<\infty$. 
    Let $(c_0^N)_{N\in\bbN}$ with $c^N_0\in\rmL^1(X_N)$ and $c^N_0\ge0$ be well-prepared, i.e., let $\iota_Nc^N_0\to \rho_0$ in $\rmL^1(X)$ and $\calE_N(c^N_0)\to \calE(\rho_0)$ as $N\to\infty$.
    
    Then, for each $N\in\bbN$ there exists a solution $c^N$ of \eqref{eq:reaction_intro} on $\bbZ^d_N$ in the sense of Definition~\ref{def:sol_disc} with initial datum $c_0^N$. 
    
    Furthermore, (up to a subsequence) we have $\iota_N c^N\to\rho$ as in Theorem~\ref{thm:Comp_Liminf}, where $\rho$ is a solution in the weak sense and in the sense of Definition~\ref{def:sol_cont} of the gradient flow equation~\eqref{eq:reaction-cross-diffusion_intro} on $\bbT^d$ with initial datum $\rho_0$.
\end{corollary}

As briefly mentioned in the introduction, it is possible to also include a drift potential into our model, thereby tilting the energy landscape. 
Namely, the following more general result holds true.
\begin{theorem}\label{thm:weighted_main-results}
Consider a \emph{reference concentration} $\omega\in\rmL^\infty(X)\cap\rmH^1(X)$ given for some potential $V\in \rmL^\infty(X)\cap\rmH^1(X)$ with $\sum_{i\in I}\exp(-V_i) = M$ a.e.\ on $\bbT^d$ by the identities $\dd\omega_i \coloneqq \exp(-V_i)\dd \scrL^d$ for all $i\in I$.
Let $\calR$, $\calR^\ast$, $\calR_N$, and $\calR^\ast_N$ as before.
Consider for the continuum model the weighted energy and slope functionals
\begin{align*}
    \cE(\rho) &\coloneqq\sum_{i\in I}\int_{\bbT^d}\LB \bra{\rho_i/\omega_i}\omega_i\dd x,\\
    \calS_\swap(\rho) &\coloneqq  \sum_{i,\ibar\in I} 2 K_{i\ibar}\int_{\bbT^d} \omega_i\omega_\ibar \abs[\big]{\sqrt{\rho_\ibar/\omega_\ibar}\nabla\sqrt{\rho_i/\omega_i}- \sqrt{\rho_i/\omega_i}\nabla\sqrt{\rho_\ibar/\omega_\ibar}}^2\dd x,\\
    \calS_\react(\rho) &\coloneqq \sum_{r\in R} 2\kappa_r\int_{\bbT^d}\omega^{(\alpha^r+\beta^r)/2}  \abs[\big]{(\rho/\omega)^{\alpha^r/2} - (\rho/\omega)^{\beta^r/2}}^2\dd x,
\end{align*}
and, setting $\discref^N \coloneqq \iota_N^\ast \omega$, for the discrete model consider the weighted functionals
\begin{align*}
    \calE_N(c) &\coloneqq \frac{1}{N^d}\sum_{i\in I}\sum_{k\in\bbZ_N^d}
    \LB \bra{c_{i,k}/\discref_{i,k}}\discref_{i,k},\\
    \calS_{N,\swap}(c) &\coloneqq  \frac{1}{N^d} \sum_{k\in\bbZ_N^d} 
    \sum_{i,\ibar\in I}  \sum_{e\in \Eu} 2 N^2 K_{i\ibar} \sqrt{\discref^N_{i,k}\discref^N_{i,k+e}\discref^N_{\ibar,k}\discref^N_{\ibar,k+e}}
    \abs[\bigg]{\sqrt{\frac{c_{i,k+e}c_{\ibar,k}}{\discref^N_{i,k+e}\discref^N_{\ibar,k}}} - \sqrt{\frac{c_{i,k}c_{\ibar,k+e}}{\discref^N_{i,k}\discref^N_{\ibar,k+e}}}}^2,\\
    \calS_{N,\react}(c) &\coloneqq  \frac{1}{N^d} \sum_{k\in\bbZ_N^d} 
    \sum_{r\in R} 2\kappa_r (\discref^N_k)^{(\alpha^r+\beta^r)/2}
    \abs[\big]{\bra{c_k/\discref_k}^{\alpha^r/2} 
    - \bra{c_k/\discref_k}^{\beta^r/2}}^2.
\end{align*}
Then, the according generalizations of Theorem~\ref{thm:Comp_Liminf}, Theorem~\ref{thm:CRContinuous}, and Corollary~\ref{cor:conv_sol} hold true.
\end{theorem}

The proofs of the main results are given in the remaining sections:
In Section \ref{s:Compactness}, we show the compactness claimed in Theorem~\ref{thm:Comp_Liminf}, while its lower limit inequalities are content of Section~\ref{s:Liminf}.
Section~\ref{s:CR} then contains the proof of Theorem~\ref{thm:CRContinuous}.
Finally, in Section~\ref{s:tilting}, we discuss the necessary adaptations of the previous arguments to obtain Theorem~\ref{thm:weighted_main-results}.

\section{Compactness}
\label{s:Compactness}

This section is dedicated to proving the first part of Theorem~\ref{thm:Comp_Liminf}, the compactness of families of curves satisfying the $N$-dependent prelimit continuity equations as well as uniform energy and dissipation bounds. 
We first embed these families into a shared space, in which we then show the claimed compactness.
A major complication that has to be overcome in Proposition~\ref{prop:strong_compactness} is the possibly vanishing mobility for species that are not fully connected.
In this case, we cannot apply an Aubin-Lions-Simon-type lemma and instead make use of the Kolmogorov-Riesz-Fréchet theorem as well as a splitting into near and far field to account for the non-quadratic nature of the prelimit structure.

Throughout this section, we fix a time horizon $T>0$.
We start our considerations from the $N$-uniform $\rmL^\infty$-bound on the energies and the $N$-uniform bound on the dissipations. We introduce explicit constants that will make it easier to track the estimates  throughout this section. 
We start with the a priori bounds
\begin{subequations}
    \label{eq:AprioriEstimBasic}
\begin{align}
\CEn &\coloneqq \sup_{N\in\bbN}\sup_{t\in[0,T]} \calE_N(c^N(t)) < \infty,
	\label{eq:AprioriEstimBEnergy}\\
\CDis &\coloneqq\sup_{N\in\bbN}\calD_N(c^N,F^N,J^N) < \infty.
	\label{eq:AprioriEstimBDiss}
\end{align}
\end{subequations}
In particular, these imply
\begin{subequations}	
\label{eq:AprioriEstimatesDissipation_parts}
\begin{align}
	\label{eq:AprioriEstimatesL1Bound}
     \CLone & \coloneqq \sup_{N\in \bbN} \sup_{t\in [0,T]} \| c^N\|_{\rmL^1_N} < \infty,
 \\   
 \CRate    &\coloneqq\sup_{N\in\bbN}\int_0^T \calR_N(c^N(t),F^N(t),J^N(t))\dd t < \infty,
	\label{eq:AprioriEstimatesRateDiff}
\\
 \CSlope	& \coloneqq \sup_{N\in\bbN}\int_0^T \calS_N(c^N(t))\dd t < \infty.
	\label{eq:AprioriEstimatesSlopeDiff}
\end{align}
\end{subequations}
Using the representation~\eqref{eq:S_swap-a-mean} and Lemma~\ref{lem:discrete_calculations_sum}, below, the slope control~\eqref{eq:AprioriEstimatesSlopeDiff} implies for $U = \sqrt{c}$ the additional uniform bounds
\begin{subequations}
\label{eq:AprioriRegularity}
\begin{align}
        \label{eq:AprioriRegularity_V}
        \CBab &\coloneqq \sup_{N\in\bbN}\max_{\ibar\in \Ifull}\int_0^T\frac{1}{N^d} \sum_{k\in\bbZ_N^d} \sum_{e\in \Eu} N^2 \abs[\big]{\babla U^N_{\ibar,k,e}(t)}^2\dd t <\infty,\\
	\label{eq:AprioriRegularity_V2}
	\CBabsq &\coloneqq \sup_{N\in\bbN}\max_{\ibar\in \Ifull}\int_0^T\frac{1}{N^d} \sum_{k\in\bbZ_N^d} \sum_{e\in \Eu} N^2 \abs[\big]{\babla c^N_{\ibar,k,e}(t)}^2\dd t <\infty,\\
	\label{eq:AprioriRegularity_U}
	\CProd &\coloneqq \sup_{N\in\bbN}\max_{\ibar\in \Ifull}\max_{i\in I} \int_0^T\frac{1}{N^d} \sum_{k\in\bbZ_N^d} \sum_{e\in \Eu}
	N^2 \abs[\big]{\tfrac12 \pra[\big]{U^N_{\ibar,k}(t)+U^N_{\ibar,k+e}(t)} \babla U^N_{i,k,e}(t)}^2\dd t <\infty.
\end{align}

\end{subequations}
 
\begin{lemma}\label{lem:discrete_calculations_sum}
        Recall the discrete gradient $\babla U_{k,e} = U_{k+e}-U_k$ and denote the arithmetic mean by $\theta(U)_{k,e} = (U_k+U_{k+e})/2$.
        Fix $\ibar\in \Ifull$ and denote $\bar U = U_\ibar$. 
        Then, it holds on $\bbZ^d_N\times\Eu$
        \begin{align*}
            \sum_{i\in I} K_{i\ibar}\abs[\big]{\theta(\bar U) \babla U_i  - \theta(U_i) \babla \bar U}^2 
            &\ge K_\ast\sum_{i\ne\ibar}\abs[\big]{\theta(\bar U) \babla U_i}^2
            + \frac{K_\ast\maxdensity}{2}\abs[\big]{\babla \bar U}^2.
        \end{align*}
    \end{lemma}
    \begin{proof}
        By definition we have $\sum_{i\ne \ibar}U^2_i + \bar U^2 = M$ everywhere on $\bbZ^d_N$. 
        This in particular implies $\sum_{i\ne\ibar}\babla (U^2_i) = -\babla (\bar U^2)$ everywhere on $\bbZ^d_N\times\Eu$. 
        Combining this with the elementary identity $\theta(U)\babla U = \frac12 \babla U^2$, we compute
\begin{align*}
        \sum_{i\in I} \abs[\big]{\theta(\bar U) \babla U_i  - \theta(U_i) \babla \bar U}^2
        &= \sum_{i\ne\ibar} \abs[\big]{\theta(\bar U) \babla U_i  - \theta(U_i) \babla \bar U}^2 \\
        &= \sum_{i\ne\ibar}\pra[\Big]{\abs[\big]{\theta(\bar U) \babla U_i}^2 
        + \abs[\big]{\theta(U_i) \babla \bar U}^2}
        + \frac12 \abs[\big]{\babla (\bar U)^2}^2\\
        &= \sum_{i\ne\ibar}\abs[\big]{\theta(\bar U) \babla U_i}^2
        + \pra[\Big]{\sum_{i\ne\ibar}\theta^2(U_i)+2\theta^2(\bar U)} \abs[\big]{\babla \bar U}^2\\
        &\ge \sum_{i\ne\ibar}\abs[\big]{\theta(\bar U) \babla U_i}^2
        + \frac{\maxdensity}{2}\abs[\big]{\babla \bar U}^2.
\end{align*}
This concludes upon recalling the lower bound in assumption~\eqref{as_general:pos_coeff_weak}.
\end{proof}

Using the embeddings from Section \ref{s:Embedding}, we define the curves 
\[
\rho^N\coloneqq\iota_Nc^N, \quad f^N\coloneqq\iota_{N,\swap}F^N, \quad j^N\coloneqq\iota_{N,\react}J^N,
\]
and show that these have weak limits.
\begin{proposition}[Weak limits]
\label{prop:FluxBound} 
Let $(c^N,F^N,J^N)$ satisfy the a priori estimates \eqref{eq:AprioriEstimBasic}. 
Then, for all $N\in\bbN$, $e\in\Eu$, $i\in I$, $\ibar\in\Ifull$, $r\in R$,
we have
\begin{equation}
\label{eq:FluxBound}
\begin{aligned}
\frac{1}{N^d}
\sum_{k\in\bbZ^d_N}
\int_0^T \C\bra[\bigg]{\frac{F^N_{i,\ibar,k,e}}{N}} \dd t 
\le \frac{1}{N^d}
\sum_{k\in\bbZ^d_N}
\int_0^T N^2\C\bra[\bigg]{\frac{F^N_{i,\ibar,k,e}}{N^2}} \dd t 
&\le C_{M,K}\CRate,\\
\frac{1}{N^d}
\sum_{k\in\bbZ^d_N}\int_0^T \C\bra{J^N_{r,k}} \dd t 
&\le C_{M,\kappa} \CRate,
\end{aligned}
\end{equation} 
where $C_{M,K}>0$ depends only on $M$, $K_*$, $K^*$ and $C_{M,\kappa}>0$ only on $M$, $\min_{r\in R}\kappa_r$.
Moreover, there exist curves of fluxes $(f,j)\in\rmL^1([0,T]\times Y)$ 
s.t.\ (along a subsequence), we have
\begin{align*}
	\int_0^T\int_{\bbT^d} \C(f_{i,\ibar,e}) \dd x \dd t \leq C_{M,K} \CRate 
	\qquad &\text{and} \qquad 
	\int_0^T\int_{\bbT^d} \C(j_r) \dd x \dd t \leq C_{M,\kappa}\CRate,\\
	 f^N_{i,\ibar,e} \rightharpoonup f_{i,\ibar,e} \quad\text{weakly in $\rmL^1([0,T]\times Y_\swap)$}
	 \qquad &\text{and} \qquad
	  j^N_r \rightharpoonup j_r
      \quad\text{weakly in $\rmL^1([0,T]\times Y_\react)$}.
\end{align*}
Finally, there exists $\rho \in \rmL^\infty([0,T]\times X)\cap \AC([0,T];\calM_+(X))$ s.t.\ (along a subsequence) it holds $\rho^N\rightharpoonup^\ast \rho$ weakly-$^\ast$ in $\rmL^\infty([0,T]\times X)$ and $\rho^N(t)\rightharpoonup^\ast\rho(t)$ weakly-$^\ast$ in $\calM_+(X)$ for all $t\in[0,T]$.
\end{proposition}
\begin{proof} 
We consider only the swap fluxes, noting that owing to the uniform $\rmL^\infty$-control over $c^N$ treating the reaction fluxes is even simpler due to the absence of the leading factor $N^2$. 
We fix indices $i\in I$, $\ibar\in\Ifull$, and $e\in\Eu$ and suppress them throughout the proof. 
Since for all $a\in\R$ the map $b\mapsto \CC(a|b)$ is monotone decreasing, the uniform $L^\infty$-bound on $c^N$ directly yields
\begin{align*}
   &\frac{1}{N^d}\sum_{k\in\bbZ^d_N}\int_0^T N^2\C\bigg(\frac{F^N_k}{N^2}\bigg) \dd t
    \le \frac{C_{M,K}}{N^d}\sum_{k\in\bbZ^d_N}\int_0^T N^2\CC\bigg(\frac{F^N_k}{N^2} \bigg| \mob_k(c^N) \bigg) \dd t\\
    &\qquad=  \frac{C_{M,K}}{N^d}\sum_{k\in\bbZ^d_N}\int_0^T \CC\big( F^N_k \big| N^2\mob_k(c^N) \big) \dd t \le \int_0^T \calR_N(c^N,F^N,J^N)\dd t \le C_{M,K}\CRate.
\end{align*} 
This shows the second inequality in the first line of \eqref{eq:FluxBound}. 
The first inequality directly follows from the fact that $b\mapsto b^2 \rmC(a/b)$ is monotone increasing.
With this, the compactness for $f$ follows from \eqref{eq:embedded-flux-bound} and the superlinear growth of $\C$.

Regarding the time regularity for $\rho^N$, we observe that Lemma~\ref{lem:CE.DiscCont} and \eqref{eq:embedded-flux-bound} imply for every $\varphi\in \rmC^1(X)$
\begin{align*}
        \langle (F^N,J^N),\nablaDisc\iota_N^\ast \varphi \rangle
        &= \langle (f^N,j^N),\nablaCont \varphi \rangle
        \le C\norm{\varphi}_{\rmW^{1,\infty}(X)} 
        \bra[\Bigg]{\sum_{i\in I}\sum_{\ibar\in \Ifull}\sum_{e\in \Eu}\norm[\bigg]{\frac{F_{i,\ibar,e}^N}{N}}_{\rmL^1_N} + \sum_{r\in R} \norm{J^N_r}_{\rmL^1_N}},
    \end{align*}
    and after an integration in time, the flux terms on the right-hand side are bounded since by convexity $s\xi \le \rmC(s) + \rmC^\ast(\xi)$ for all $s,\xi \ge0$ (in particular for $\xi = 1$) and using that $\rmC^\ast(1) \le 1$.
		
    Now we consider any partition $(t_l)_{l=0}^{l^\ast}$, $l^\ast\in\bbN$ of $[0,T]$. 
    Then, the previous bounds and the discrete continuity equation yield for every $\varphi\in \rmC^1(X)$ the estimate 
    \begin{align*}
         \sum_{l=1}^{l^\ast} \langle \rho^N(t_l)
           -\rho^N(t_{l-1}), \varphi\rangle
        &= \sum_{l=1}^{l^\ast} \langle c^N(t_l)-c^N(t_{l-1}),
              \iota_N^\ast\varphi\rangle_N
        \le C\norm{\varphi}_{\rmW^{1,\infty}(X)}.
    \end{align*}
Taking suprema with respect to $\varphi\in\set{\varphi\in \rm\rmC^1(X)$: $\norm{\varphi}_{\rmW^{1,\infty}(X)}\le 1}$ and the partition, we obtain $\sup_N\norm{\rho^N}_{\mathrm{BV}([0,T];(\mathrm{W}^{1,\infty}(X))^\ast)}<\infty$ and an application of \cite[Proposition~3.3.1]{AGS} yields
$\rho\in\AC([0,T];\calM_+(X))$ s.t.\ (up to a subsequence) $\rho^N(t) \rightharpoonup^\ast \rho(t)$ in the weak-$^\ast$ sense of measures for all $t\in[0,T]$, while weak-$^\ast$ convergence in $\rmL^\infty([0,T]\times X)$ follows (after taking another subsequence) from the uniform $\rmL^\infty$-boundedness of $\rho^N$.
\end{proof}

The above modes of convergence immediately imply that any limit is in $\CE^\maxdensity$.
\begin{corollary}[Closedness of $\CE^\maxdensity$]\label{cor:CE_closed}
    Let $(\rho,f,j)$ be a limit of $(\rho^N,f^N,j^N)_{N\in\bbN}\subset\CE^\maxdensity$ in the sense of Proposition~\ref{prop:FluxBound}. Then, it holds $(\rho,f,j)\in\CE^\maxdensity$. 
\end{corollary}
\begin{proof}
    By definition the set $\CE^\maxdensity$ is closed with respect to the weak $\rmL^1$ convergence of the time-integrated embedded fluxes and the pointwise-in-time weak-$^\ast$ convergence of the embedded concentrations shown in Proposition~\ref{prop:FluxBound}.
\end{proof}
    
Our next goal is establishing spatial regularity of the limit curves. 
We begin with the fully connected species $\rho^N_\idot$, $i\in\Ifull$ by building on \eqref{eq:AprioriRegularity}.
\begin{proposition}[Spatial regularity and strong convergence for fully connected $\rho_i$] \label{prop:space_reg_bar-rho}
    Let $(c^N)_N$ satisfy \eqref{eq:AprioriEstimBasic}, $\rho$ be the limit of $\rho^N$ as in Proposition~\ref{prop:FluxBound}, and $i\in \Ifull$.
    Then, the embedded sequence $(\rho^N_\idot)_N$ with $\rho^N_\idot = \iota_N c^N_\idot$ is uniformly bounded in $\rmL^1([0,T];\rmB\rmV(X))$
    and we have the convergence $\rho^N_\idot \to \rho_i$ strongly in $\rmL^p([0,T]\times\bbT^d)$ for all $p\in[1,\infty)$.
    Furthermore, the limit has the regularity $\sqrt{\rho_i}\in\rmL^2([0,T];\rmH^1(X))$.
\end{proposition}
\begin{proof}
	We fix and drop the index $i \in\Ifull$ throughout the proof. 
    Due to the uniform $\rmL^\infty$-bound it suffices to consider $p =1$.
	We begin with the following observation: For $x,y\in\bbT^d$ define the set of cell boundaries intersecting the line between $x-y$ and $x$
    \begin{align*}
        \mathrm{Int}(x,y) \coloneqq \set{(k,e)\in\bbZ^d_N\times\Eu: [x-y,x]\cap\partial Q^N_{k+e}\cap\partial Q^N_k \ne \emptyset}.
    \end{align*}
    On the other hand, for $y\in\bbT^d$, $k\in\bbZ^d_N$, and $e\in\Eu$, we define the set of points $x$ such that the line between $x-y$ and $x$ intersects the boundaries between $Q^n_k$ and $Q^N_{k+e}$ 
    \begin{align*}
        \mathrm{Int}(y,k,e) \coloneqq \set{x\in\bbT^d: [x-y,x]\cap\partial Q^N_{k+e}\cap\partial Q^N_k \ne \emptyset}. 
    \end{align*} 
    Then, it holds for all $x,y\in\bbT^d$ 
    \begin{align}
    \label{eq:cylinder_estimate}
        \abs{\rho^N(x-y)-\rho^N(x)} 
        &\le \sum_{k\in\bbZ^d_N}\sum_{e\in\Eu}\abs{c^N_{k+e}-c^N_k}\one_{\mathrm{Int}(x,y)}(k,e)
        =\sum_{k\in\bbZ^d_N}\sum_{e\in\Eu}\abs{c^N_{k+e}-c^N_k}\one_{\mathrm{Int}(y,k,e)}(x).
    \end{align}
    Noticing that $\mathrm{Int}(y,k,e)$ is an (oblique) cylinder with base $\partial Q^N_{k+e}\cap\partial Q^N_k$ (which has area $1/N^{d-1}$) and axis parallel to $y$, integrating this inequality w.r.t. $x$, and summing over all pairs of cells, we obtain for any $y\in\bbT^d$
	 \begin{align*}
        \int_0^T\int_{\bbT^d}\abs{\rho^N(x-y)-\rho^N(x)}\dd x\dd t
        &\le \frac{\abs{y}}{N^{d-1}} \int_0^T\sum_{k\in\bbZ^d_N}\sum_{e\in\Eu}\abs[\big]{c^N_{k+e}-c^N_{k}}\dd t\\
        &= \frac{\abs{y}}{N^d} \int_0^T\sum_{k\in\bbZ^d_N}\sum_{e\in\Eu} N\abs[\big]{U^N_{k+e}-U^N_k}\bra[\big]{U^N_{k+e} + U^N_k} \dd t\\
        &\le 2\sqrt{\CBab T}\abs{y}.
    \end{align*}    
    In particular, for every $\varphi\in C^1([0,T]\times\bbT^d)$ and every $l = 1,\ldots, d$ we can approximate $\partial_{x_l}\varphi$ by difference quotients to obtain
    \begin{align*}
        \int_0^T\int_{\bbT^d}\rho^N \partial_{x_l}\varphi \dd x\dd t 
        \le  2\sqrt{\CBab T}\norm{\varphi}_{L^\infty},
    \end{align*}
    i.e., $\rho^N$ is uniformly bounded in $\rmL^1([0,T];\rmB\rmV(\bbT^d))$.

    Combining this with the time regularity established in Proposition~\ref{prop:FluxBound}, the $\rmL^1$-compactness for $\rho^N$ follows as an application of \cite[Theorem~2 and Proposition~1.10]{rossi2003tightness} by arguing exactly as in the proof of \cite[Theorem~4.8]{HT23}.

    From this, the $\rmL^2$-compactness of $\sqrt{\rho^N}$ is an immediate consequence. Indeed, since $\rho^N,\rho$ are non-negative, we have $\abs{\sqrt{\rho^N}-\sqrt{\rho}}^2 \le \abs{\rho^N-\rho}$.
    
    Finally, the spatial regularity of $\sqrt{\rho}$ for any limit curve $\rho$ follows from the uniform control over the finite differences in terms of the slope, where the strong $\rmL^2$-convergence $\sqrt{\rho^N} \to \sqrt{\rho}$ allows us to identify the limit as $\nabla\sqrt{\rho}$. (cf.\ \cite[Proposition~5.11]{HMS25} for more details). 
\end{proof}

Next, we would like to derive strong compactness for species $i\in I\setminus\Ifull$.
However, the a priori bound \eqref{eq:AprioriRegularity_U} contains a (possibly vanishing) mobility. 
Therefore, a classical Aubin-Lions-Simon-type approach becomes unfeasible
and we will instead argue based on the Kolmogorov-Riesz-Fréchet theorem similar to \cite[following (4.33)]{burger2010nonlinear}.
However, since we consider not only gradient flow solutions, but all curves with bounded dissipation and energy, we cannot rely on quadratic estimates such as the Cauchy-Schwarz inequality.
Indeed, while it does hold $\C(r) \approx \frac{1}{2}r^2$ near zero, for $r\gg 1$ we only have $\C(r)\approx 2 r\log r$.
In the proof we will thus treat these two regimes differently.

\begin{proposition}[Weighted strong convergence for partially connected species]\label{prop:strong_compactness}
    Let $(c^N)_N$ satisfy \eqref{eq:AprioriEstimBasic}. 
    Let $\rho$ be the limit of $\rho^N$ as in Proposition~\ref{prop:FluxBound}.
    Then, for any $i\in I$, $\ibar\in\Ifull$, $p\in[1,\infty)$, 
    we have the strong convergences
    \begin{subequations}
    \begin{align}
        \label{eq:strong_limit_product}
        \rho^N_\ibar\rho^N_i &\to \rho_\ibar\rho_i 
        \qquad\; \text{in $\rmL^p([0,T]\times \bbT^d)$},\\
        \label{eq:strong_limit_indicator}
        \rho^N_\idot\one_{\set{\rho_\ibar>0}} &\to \rho_i\one_{\set{\rho_\ibar>0}} \quad \text{in $\rmL^p([0,T]\times \bbT^d)$}.
    \end{align}
    \end{subequations}
\end{proposition}
\begin{proof}
    Throughout the proof we fix and suppress an index $i\in I$. 
    W.l.o.g.\ we assume that $K_{i0} > 0$ and $K_{i\ibar} = 0$ for $i,\ibar\in I\setminus\set{0}$ and we denote $\bar\rho = \rho_0$.
    Indeed, adding more fully connected species does not change any of the arguments, but makes the notation more cumbersome.

    We first consider~\eqref{eq:strong_limit_product} and due to the $\rmL^\infty$-bound w.l.o.g.\ only show the strong $\rmL^2$-convergence $\sqrt{\rho^N_\ibar}\rho^N_i \to \sqrt{\rho_\ibar}\rho_i$. 
    By the Kolmogorov-Riesz-Fréchet theorem (cf.\ e.g.\ \cite[Theorem~4.26]{brezis2011functional}) and the uniform bounds in $\rmL^\infty$, this further reduces to showing
    \begin{align*}
        \lim_{\tau,\abs{y}\to 0}\norm[\big]{\bra[\big]{\sqrt{\bar\rho^N}\rho^N}(\cdot+\tau,\cdot+y) - \sqrt{\bar\rho^N}\rho^N}_{\rmL^2([0,T]\times\bbT^d)} = 0.
    \end{align*}
    To this end, we fix $y\in\R^d$, $\tau >0$ and (after setting $\rho^N(t) = \rho^N(2T-t)$ for $t\in [T,T+\tau]$) we estimate
    \begin{align*}
        &\int_0^T \int_{\bbT^d}\abs[\big]{\sqrt{\bar\rho^N}\rho^N(t+\tau,x+y) - \sqrt{\bar\rho^N}\rho^N(t,x)}^2 \dd x \dd t\\
        &\le 2 \underbrace{\int_0^T\int_{\bbT^d}\abs[\big]{\sqrt{\bar\rho^N}\rho^N(t+\tau,x+y) - \sqrt{\bar\rho^N}\rho^N(t+\tau,x)}^2 \dd x \dd t}_{\eqqcolon \rmI}\\
        &+ 2 \underbrace{\int_0^T\int_{\bbT^d}\abs[\big]{\sqrt{\bar\rho^N}\rho^N(t+\tau,x) - \sqrt{\bar\rho^N}\rho^N(t,x)}^2 \dd x \dd t}_{\eqqcolon \rmI\rmI}.
    \end{align*}
    For the first term, we estimate using \eqref{eq:cylinder_estimate} and \eqref{eq:AprioriRegularity} to derive
    \begin{align*}      
        \rmI &\le \int_0^T\frac{\abs{y}}{N^{d-1}}\sum_{k\in\bbZ^d_N}\sum_{e\in\Eu} \abs[\Big]{\sqrt{\bar c^N_{k+e}} c^N_{k+e} - \sqrt{\bar c^N_k} c^N_k}^2(t+\tau) \dd t\\
        &\le \int_0^T\frac{\abs{y}}{N^{d-1}}\sum_{k\in\bbZ^d_N}\sum_{e\in\Eu} 2\maxdensity\bra[\Big]{\maxdensity\abs[\Big]{\sqrt{\bar c^N_{k+e}} - \sqrt{\bar c^N_k}}^2 
        + 2\bar c^N_k 
        \abs[\Big]{
        \sqrt{c^N_{k+e}} - \sqrt{c^N_k}
        }^2}(t+\tau) \dd t\\
        &\le C
        \abs{y} \overset{y\to 0}{\longrightarrow}0.
    \end{align*}
    For the second term we split the expression again, writing
    \begin{align*}
        \rmI\rmI &\le 2 \underbrace{\int_0^T\int_{\bbT^d}\bra[\big]{\rho^N}^2(t+\tau,x)\abs[\big]{\sqrt{\bar\rho^N}(t+\tau,x) - \sqrt{\bar\rho^N}(t,x)}^2 \dd x \dd t}_{\eqqcolon \rmI\rmI_1}\\
        &+ 2 \underbrace{\int_0^T\int_{\bbT^d}\bar\rho^N(t,x)\abs[\big]{\rho^N(t+\tau,x) - \rho^N(t,x)}^2 \dd x \dd t}_{\eqqcolon \rmI\rmI_2}.
    \end{align*}
    For $\rmI\rmI_1$, we use $\rmL^\infty$-bound on $\rho^N$ to estimate
    \begin{align*}
        \rmI\rmI_1 \le M^2\int_0^T\int_{\bbT^d}\abs[\big]{\sqrt{\bar\rho^N}(t+\tau,x) - \sqrt{\bar\rho^N}(t,x)}^2 \dd x \dd t.
    \end{align*}
    Since $\sqrt{\bar\rho^N}$ is compact $\rmL^2([0,T]\times\bbT^d)$ as shown in Proposition~\ref{prop:space_reg_bar-rho}, an application of the Kolmogorov-Riesz-Fréchet theorem (in reversed direction) implies that the upper bound vanishes uniformly in $N$ as $\tau\to 0$.
    
    For $\rmI\rmI_2$ we aim to
    argue using the continuity equation, in particular building on \eqref{eq:FluxBound}. To this end, we split into linear differences
    \begin{align*}
        \rmI\rmI_2 &= - \underbrace{\int_0^T\int_{\bbT^d}\bar\rho^N(t,x)\rho^N(t,x)\bra[\big]{\rho^N(t+\tau,x) - \rho^N(t,x)} \dd x \dd t}_{\eqqcolon \rmI\rmI_{21}}\\
        &+\underbrace{\int_0^T\int_{\bbT^d}\bar\rho^N(t,x)\rho^N(t+\tau,x)\bra[\big]{\rho^N(t+\tau,x) - \rho^N(t,x)} \dd x \dd t}_{\eqqcolon \rmI\rmI_{22}}
    \end{align*}
    For $\rmI\rmI_{21}$ we apply the discrete continuity equation and split the arising swap terms once more into sets depending on a parameter $\alpha>0$ to be determined later 
    \begin{align*}
        \rmI\rmI_{21} &= \int_0^T\frac{1}{N^d}\sum_{k\in\bbZ^d_N}\bar c^N_k(t)c^N_k(t)(c^N_k(t+\tau)-c^N_k(t)) \dd t\\
        &= \int_0^T\int_t^{t+\tau}\frac{1}{N^d}\sum_{k\in\bbZ^d_N}\bra[\Bigg]{\sum_{e\in\Eu}\nablaDDiff[\bar c^N(t)c^N(t)]_{k,e} F^N_{k,e}(s) + \sum_{r\in R}\nablaDReact[\bar c^N(t)c^N(t)]_{r,k}J_{r,k}(s) } \dd s \dd t\\
        &= \underbrace{\int_0^T\int_t^{t+\tau}\frac{1}{N^d}\sum_{k\in\bbZ^d_N}\sum_{e\in\Eu}\nablaDDiff[\bar c^N(t)c^N(t)]_{k,e} F^N_{k,e}(s)\one_{\set{\abs{F^N_{k,e}} < N^2\tau^{-\alpha}}}(s) \dd s \dd t}_{\eqqcolon \rmI\rmI_{211}}\\
        &+\underbrace{\int_0^T\int_t^{t+\tau}\frac{1}{N^d}\sum_{k\in\bbZ^d_N}\sum_{e\in\Eu}\nablaDDiff[\bar c^N(t)c^N(t)]_{k,e} F^N_{k,e}(s)\one_{\set{\abs{F^N_{k,e}} \ge N^2\tau^{-\alpha}}}(s) \dd s \dd t}_{\eqqcolon \rmI\rmI_{212}}\\
        &+\underbrace{\int_0^T\int_t^{t+\tau}\frac{1}{N^d}\sum_{k\in\bbZ^d_N}\sum_{r\in R}\nablaDReact[\bar c^N(t)c^N(t)]_{r,k}J_{r,k}(s)\dd s \dd t}_{\eqqcolon \rmI\rmI_{213}}.
    \end{align*}
    The term $\rmI\rmI_{213}$ converges uniformly to zero due to the uniform $\rmL^\infty$-bound on $c^N$ and the uniform integrability of $J^N$ established in Proposition~\ref{prop:FluxBound}.
    
    To control $\rmI\rmI_{211}$ we recall that by~\eqref{eq:AprioriRegularity} the term $N \nablaDDiff[\bar c^Nc^N]_{k,e} = N \babla\bra{\bar c^Nc^N}_{k,e} - N \babla\bra{\bar c^N}^2_{k,e}$
	is uniformly bounded in $\rmL^2$.
	Therefore, since $\C(r) \approx r^2$ for $0\le r\le 1$ and since the map $s\mapsto s\C(r/s)$ is decreasing, Cauchy-Schwarz and Proposition~\ref{prop:FluxBound} yield for any $\alpha,\tau\in(0,1)$
    \begin{align*}
        \rmI\rmI_{211}
        &\le \bra[\Bigg]{\tau\int_0^T\frac{1}{N^d}\sum_{k\in\bbZ^d_N}\sum_{e\in\Eu} N^2\abs[\big]{\nablaDDiff[\bar c^N(t)c^N(t)]_{k,e}}^2\dd t}^{\frac{1}{2}}\\
        &\qquad\qquad\quad\cdot
        \bra[\Bigg]{T\int_0^{T+\tau}\frac{1}{N^d}\sum_{k\in\bbZ^d_N}\sum_{e\in\Eu}N^2\tau^{-2\alpha}\abs[\bigg]{\frac{F^N_{k,e}(s)}{N^2\tau^{-\alpha}}}^2\one_{\set{\abs{F^N_{k,e}} < N^2\tau^{-\alpha}}}(s) \dd s}^{\frac{1}{2}} 
        \displaybreak[0]
        \\
        &\le 
        \sqrt{C_2 T}\tau^{\frac{1-\alpha}{2}}
        \bra[\Bigg]{\int_0^{T+\tau}\frac{1}{N^d}\sum_{k\in\bbZ^d_N}\sum_{e\in\Eu} N^2\tau^{-\alpha}\C\bra[\bigg]{\frac{F^N_{k,e}(s)}{N^2\tau^{-\alpha}}} \one_{\set{\abs{F^N_{k,e}} < N^2\tau^{-\alpha}}}(s) \dd s}^{\frac{1}{2}}
        \displaybreak[0]
        \\
        &\le 
        \sqrt{C_3 T}\tau^{\frac{1-\alpha}{2}}
        \bra[\Bigg]{\int_0^{T+\tau}\frac{1}{N^d}\sum_{k\in\bbZ^d_N}\sum_{e\in\Eu} N^2\C\bra[\bigg]{\frac{F^N_{k,e}(s)}{N^2}} \dd s}^{\frac{1}{2}}
        \le 
        \sqrt{C T}\tau^{\frac{1-\alpha}{2}}
        \overset{\tau \to 0}{\longrightarrow}0,
    \end{align*}
    where we used that $F^N_{k,e}(t) = -F^N_{k,e}(2T-t)$ for $t\in[T,T+\tau]$ by the choice of extension for $\rho^N$.
    
    Next, we consider $\rmI\rmI_{212}$. Here, arguing directly via uniform integrability is not enough due to the slow growth $\C(r)\approx r(\log r+1)$ for $r>1$ in combination with the fact that $F^N$ is not uniformly bounded without a factor $1/N$.
    To accommodate for this, we have introduced the factor $\tau^\alpha$.
    Indeed, we use the $\rmL^\infty$-control over $c^N$, that the map $r\mapsto\log r$ is increasing and positive on $(1,\infty)$, as well as the sharper bound in Proposition~\ref{prop:FluxBound} to derive for $\alpha,\tau \in(0,1)$
    \begin{align*}
        \rmI\rmI_{212}
        &\le 4\maxdensity T\int_0^{T+\tau}\frac{1}{N^d}\sum_{k\in\bbZ^d_N}\sum_{e\in\Eu} N^2\frac{\abs[\big]{F^N_{k,e}(s)}}{N^2}\one_{\set{\abs{F^N_{k,e}}/N^2 \ge \tau^{-\alpha}}}(s) \dd s\\
        &\le \frac{4\maxdensity T}{-\alpha\log\tau} \int_0^{T+\tau}\frac{1}{N^d}\sum_{k\in\bbZ^d_N}\sum_{e\in\Eu} N^2\frac{\abs[\big]{F^N_{k,e}(s)}}{N^2}\log\bra[\bigg]{\frac{\abs[\big]{F^N_{k,e}(s)}}{N^2}}\one_{\set{\abs{F^N_{k,e}}/N^2 \ge \tau^{-\alpha}}}(s) \dd s\\
        &\le \frac{4\tilde C\maxdensity T}{\alpha\abs{\log\tau}} \int_0^{T+\tau}\frac{1}{N^d}\sum_{k\in\bbZ^d_N}\sum_{e\in\Eu} N^2\C\bra[\bigg]{\frac{F^N_{k,e}(s)}{N^2}} \dd s
        \le \frac{4C\maxdensity T}{\alpha\abs{\log\tau}} \overset{\tau\to 0}{\longrightarrow}0.
    \end{align*}
    Hence, choosing any $\alpha\in(0,1)$, both $\rmI\rmI_{211}$ and $\rmI\rmI_{212}$ vanish as $\tau\to 0$.
    
    It remains to control $\rmI\rmI_{22}$. We split this term yet again as follows:
    \begin{align*}
        \rmI\rmI_{22} 
        &= \underbrace{\int_0^T\int_{\bbT^d}\bar\rho^N(t+\tau,x)\rho^N(t+\tau,x)\bra[\big]{\rho^N(t+\tau,x) - \rho^N(t,x)}}_{\eqqcolon \rmI\rmI_{221}} \dd x \dd t\\
        &+ \int_0^T\int_{\bbT^d}\underbrace{\bra[\big]{\bar\rho^N(t,x)- \bar\rho^N(t+\tau,x)}\rho^N(t+\tau,x)\bra[\big]{\rho^N(t+\tau,x) - \rho^N(t,x)}}_{\eqqcolon \rmI\rmI_{222}} \dd x \dd t.
    \end{align*}
    $\rmI\rmI_{221}$ vanishes by the same arguments used for $\rmI\rmI_{21}$, while $\rmI\rmI_{222}$ can be treated as $\rmI\rmI_1$, since by Proposition~\ref{prop:space_reg_bar-rho} we also have strong compactness of $\bar\rho^N$ in $\rmL^1([0,T]\times\bbT^d)$. 

    Hence, by \cite[Theorem 4.26]{brezis2011functional} there exists $v\in \rmL^2([0,T]\times \bbT^d)$ s.t.\ (along a not renamed subsequence) $\sqrt{\bar\rho^N}\rho^N \to v$ in $\rmL^2([0,T]\times \bbT^d)$. Proposition~\ref{prop:FluxBound} and Proposition~\ref{prop:space_reg_bar-rho} then allow us to identify $v = \sqrt{\bar\rho}\rho$.

    To obtain~\eqref{eq:strong_limit_indicator}, we observe that combining~\eqref{eq:strong_limit_product} with the strong $\rmL^p$-limit $\bar\rho^N\to\bar\rho$ from Proposition~\ref{prop:space_reg_bar-rho}, we have a.e.\ on $\set{\bar\rho>0}$ the pointwise limit
    \begin{align*}
        \rho^N = \frac{1}{\sqrt{\bar\rho}}\bra[\big]{\sqrt{\bar\rho}-\sqrt{\bar\rho^N}}\rho^N + \frac{1}{\sqrt{\bar\rho}}\sqrt{\bar\rho^N}\rho^N \to \rho.
    \end{align*}
    Then, the strong $\rmL^p$-limit follows via dominated convergence.
    This concludes the proof.
\end{proof}

Next, we study the spatial regularity of the individual species.
\begin{proposition}[Spatial regularity of $\rho$]\label{prop:space_reg_rho}
    Let $(c^N)_N$ satisfy \eqref{eq:AprioriEstimBasic} and let $\rho$ be the limit as in Propositions~\ref{prop:FluxBound}, \ref{prop:space_reg_bar-rho}, and \ref{prop:strong_compactness}.
    Denote $\babla_N: L^1(\bbT^d)\to L^1(\bbT^d;\bbR^d)$, $\babla_N\rho \coloneqq N\sum_{e\in\Eu} (\shift^N_e\rho - \rho)e$.
    Then, for all $i\in I$,$\ibar\in\Ifull$ we have (along a subsequence)  
    \begin{subequations}\label{eq:babla_conv}
    \begin{align}
        \label{eq:babla_conv1}
        \sqrt{\rho^N_\idot}\babla_N\sqrt{\rho^N_\ibar} \one_{\set{\rho_\ibar>0}} \rightharpoonup \sqrt{\rho_i}\nabla\sqrt{\rho_\ibar}\quad \text{weakly in }\rmL^2([0,T]\times \bbT^d;\bbR^d),\\
        \label{eq:babla_conv3}
        \sqrt{\rho^N_\ibar}\babla_N\sqrt{\rho^N_\idot} \one_{\set{\rho_\ibar>0}} \rightharpoonup \sqrt{\rho_\ibar}\nabla\sqrt{\rho_i}\quad \text{weakly in }\rmL^2([0,T]\times \bbT^d;\bbR^d).
    \end{align}
    \end{subequations}
    In particular, $\sqrt{\rho_i}\in \rmL^2([0,T];\rmH^1(\bbT^d,\rho_\ibar))$.
\end{proposition}
\begin{proof}
    Throughout the proof we fix $i\in I$, $\ibar\in\Ifull$ and denote $\rho =\rho_i$, $\bar\rho =\rho_\ibar$. 
    We will use multiple times the strategy of bounding a sequence in $\rmL^2$, so that it has a weak $\rmL^2$-limit, and then identifying its weak $\rmL^1$-limit, which coincides with the weak $\rmL^2$-limit.
    
    Regarding~\eqref{eq:babla_conv1}, note that by the uniform bound~\eqref{eq:AprioriRegularity_V} to obtain weakly converging subsequences and the strong convergence from Proposition~\ref{prop:space_reg_bar-rho} to identify the limit, we have $\babla_N\sqrt{\bar\rho^N} \rightharpoonup \nabla\sqrt{\bar\rho}$ weakly in $\rmL^2([0,T]\times\bbT^d;\R^d)$. 
    Then, using the same type of argument, the claim follows from the uniform bound~\eqref{eq:AprioriRegularity_U} and the strong convergence~\eqref{eq:strong_limit_indicator}.
    
    To obtain \eqref{eq:babla_conv3}, we again use the bound~\eqref{eq:AprioriRegularity_V} as well as the strong $\rmL^2$-convergence $\sqrt{\bar\rho^N\rho^N}\to \sqrt{\bar\rho\rho}$ to obtain the weak $\rmL^2$-convergence $\babla_N\sqrt{\rho^N\bar\rho^N} \rightharpoonup \nabla\sqrt{\rho\bar\rho}$.
    Then, we rewrite
    \begin{align*}
         \sqrt{\bar\rho^N}\babla_N\sqrt{\rho^N}\one_{\set{\bar\rho>0}} 
         &= \babla_N\sqrt{\rho^N\bar\rho^N}\one_{\set{\bar\rho>0}} 
         - \bar\shift^N\bra[\big]{\sqrt{\rho^N}}\babla_N\sqrt{\bar\rho^N}\one_{\set{\bar\rho>0}},
    \end{align*}
    where we denoted by $\bar\shift^N$ the diagonal shift matrix with entries $\bar\shift^N_{ll} = \shift^N_{e_l}$, $l=1,\ldots, d$ defined in~\eqref{eq:shift}.
    Since also $\bar\shift^N\bra[\big]{\sqrt{\rho^N}}\one_{\set{\bar\rho>0}} \to \sqrt{\rho}\one_{\set{\bar\rho>0}}$ strongly in $\rmL^2([0,T]\times\bbT^d)$, we can identify the limit of the right-hand side as before, thus proving the claim.
\end{proof}

\section{Lower limit inequalities}
\label{s:Liminf}

In the previous section we have obtained convergence of the embedded prelimit family to a limit family $(\iota_Nc^N,\iota_{N,\swap}F^N,\iota_{N,\react}J^N) \to (\rho,f,j)$ as $N\to\infty$. 
To complete the proof of Theorem~\ref{thm:Comp_Liminf}, the next step is to prove rigorous lower limit inequalities for the obtained notion of convergence. 
To this end, we prove lower limit inequalities for the rate and slope terms separately. 
The latter is an immediate consequence of Proposition~\ref{prop:space_reg_rho} as we see next.

\begin{proposition}[Liminf for slopes]\label{prop:liminf_S}
    Let $\rho^N = \iota_Nc^N \to \rho$ strongly in $\rmL^2([0,T]\times X)$ and satisfy the weak gradient convergence~\eqref{eq:babla_conv}.
    Then, it holds
    \begin{align*}
        \liminf_{N\to\infty} \calS_N(c^N) \ge \calS(\rho).
    \end{align*}
\end{proposition}
\begin{proof}
    Regarding $\calS_\swap$, recall the reformulation~\eqref{eq:S_swap-a-mean} of $\calS_{N,\swap}$, the discrete gradient vector $\babla_N\rho = N\sum_{e\in\Eu} (\shift^N_e\rho - \rho)e$, and the diagonal shift matrix $\bar\shift^N$ with $\bar\shift^N_{ll} = \shift^N_{e_l}$ for $l=1,\ldots, d$.
    With this, restricting to the smaller sets $\set{\rho_\ibar>0}$ and using the symmetry $K_{i\ibar}=K_{\ibar i}$, we obtain
    \begin{align*}
        \calS_{N,\swap}(c^N) 
        &\ge \sum_{i,\ibar\in \Ifull}
        \frac{K_{i\ibar}}{2}\int_{\bbT^d}\abs[\Big]{\bra[\Big]{\bar S^N\bra[\big]{\sqrt{\rho^N_\ibar}} {+} \mathrm{Id}\sqrt{\rho^N_\ibar}} \babla_N\sqrt{\rho^N_\idot}
        {-} \bra[\Big]{\bar S^N\bra[\big]{\sqrt{\rho^N_\idot}} {+} \mathrm{Id}\sqrt{\rho^N_\idot}} \babla_N\sqrt{\rho^N_\ibar}}^2\dd x\\
        &\quad+\sum_{i\in I}\sum_{\ibar\in \Ifull} K_{i\ibar}\int_{\bbT^d}\abs[\Big]{\bra[\Big]{\bar S^N\bra[\big]{\sqrt{\rho^N_\ibar}} {+} \mathrm{Id}\sqrt{\rho^N_\ibar}} \babla_N\sqrt{\rho^N_\idot}\one_{\set{\rho_\ibar>0}}\\
        &\qquad\qquad\qquad\qquad\qquad
        {-} \bra[\Big]{\bar S^N\bra[\big]{\sqrt{\rho^N_\idot}} {+} \mathrm{Id}\sqrt{\rho^N_\idot}} \babla_N\sqrt{\rho^N_\ibar}\one_{\set{\rho_\ibar>0}}}^2\dd x.
    \end{align*}
    Then, since the $\rmL^2$-norm is weakly lsc., the convergence~\eqref{eq:babla_conv} implies that, in the lower limit, the right-hand side bounds $\calS_{N,\swap}(c^N)$ from above.
    Indeed, the limit integrands are only supported on the sets $\set{\rho_\ibar>0}$, so that the are them invariant under multiplication by $\one_{\set{\rho_\ibar>0}}$.
    
    Regarding $\calS_\react$, we recall that by the assumption~\eqref{as_general:react} we have $R = \Rone\cup \Rtwo$.
    For $r\in \Rone$, interpolating between the strong convergence from Proposition~\ref{prop:space_reg_rho} and the uniform $\rmL^\infty$-bound, we obtain $(\rho^N)^\gamma \to \rho^\gamma$ strongly in $\rmL^1([0,T]\times\bbT^d)$ for the relevant exponents $\gamma\in\set{\alpha^r,\beta^r,(\alpha^r+\beta^r)/2}$, whence the lower limit (in fact we have the full limit for these terms).
    
    For $r\in \Rtwo$, the integrand $\rho\mapsto \abs[\big]{\rho^{\alpha^r/2} - \rho^{\beta^r/2}}^2$ is (jointly) convex and hence the integral functional is weakly lsc. 
    Indeed, for each of these terms there exist indices $i,\ibar\in I$ such that $\abs[\big]{\rho^{\alpha^r/2} - \rho^{\beta^r/2}}^2 = \abs[\big]{\rho_i^{1/2} - \rho_\ibar^{1/2}}^2$.
    This concludes the proof.    
\end{proof}

Since its structure does not change in the limit, the liminf inequality for $\calR_\react$ may be obtained using Ioffe's lemma.
To obtain the liminf inequality for the swap rates, where the structure does change, we first prove a limsup estimate for their dual dissipation potentials with sufficiently smooth test functions, in the following lemma, and then conclude using duality. 
We highlight that this argument relies only on the convex duality of $\C$ and $\C^\ast$, but not the joint convexity of $\calR_\swap$ and $\calR^\ast_\swap$, so that the argument is independent of the (non-)concavity of the strongly converging mobility.
\begin{lemma}\label{lem:EstimateDualDissipationPotential}
Let $(c^N)_{N\in\N}$ be any sequence s.t.\ $\rho^N_\idot\rho^N_\ibar \to \rho_i\rho_\ibar$ strongly in $\rmL^1(\bbT^d)$ for all $i\in I$, $\ibar\in\Ifull$. Let $\varphi\in \rmC^1(X)$. Then, it holds
    \[
    \limsup_{N\to\infty} \calR^\ast_{N,\swap}(c^N,\nablaDDiff\iota^*_N\varphi)\leq \calR_\swap^*(\rho,\nablaCDiff\varphi).
    \]
\end{lemma}
\begin{proof}
    For any $\varphi\in \rmC^1(X)$, $i\in I$, $k\in\bbZ^d_N$, and $e\in\Eu$, it holds
    \begin{align*}
        \iota^*_N\varphi_{i,k+e} - \iota^*_N\varphi_{i,k} &=
        N^d\int_{Q^N_k} \varphi_i(x-e/N) -\varphi_i(x)\dd x 
        = -\frac1N\iota^*_N(\partial_{x_e}\varphi_i)_k + o(1/N),
    \end{align*}
where the $o(1/N)$ comes from the modulus of continuity of $\partial_{x_e}\varphi$.
Hence, since $K_{i\ibar}=K_{\ibar i}$, $\C^*$ is symmetric, and $[0,\infty)\ni r\mapsto \C^*(r)$ is increasing, it follows
\begin{align*}
    \calR^\ast_{N,\swap} (c^N,\nablaDDiff\iota^*_N\varphi) &= \frac{2}{N^d} \sum_{i\in I}\sum_{\ibar\in\Ifull} \sum_{k\in\bbZ_N^d} \sum_{e\in \Eu} N^2 \mob_{i,\ibar,k,e}(c)
    \C^*(\iota^*_N\varphi_{i,k+e} {-} \iota^*_N\varphi_{i,k} {-} \iota^*_N\varphi_{\ibar,k+e} {+} \iota^*_N\varphi_{\ibar,k})\\
    &\le \frac{2}{N^d} \sum_{i\in I} \sum_{\ibar\in\Ifull} \sum_{k\in\bbZ_N^d} \sum_{e\in \Eu} N^2 \mob_{i,\ibar,k,e}(c)    \C^*\bra[\Big]{\frac1N\iota^*_N(\partial_{x_e}\varphi_i - \partial_{x_e}\varphi_\ibar)_k+o(1/N)}.
\end{align*}
Recall that, by the definition of $\C^\ast$, for all $r\in[0,\infty)$ and all $N\in\bbN$ it holds 
\[
    N^2\C^*\bra[\Big]{\frac{r}{N}}\leq \frac{r^2}{2}\cosh\bra[\Big]{\frac{r}{N}}.
\]
We want to apply this with $r=\iota^*_N(\partial_{x_e}\varphi_i - \partial_{x_e}\varphi_\ibar)_k + o(1)$, which is why we introduce the scalar
\[
    a_N\coloneqq\max_{i\in I,\ibar\in\Ifull,k\in\bbZ_N^d,e\in\Eu}\Big\{\cosh\Big(\frac1N\iota^*_N(\partial_{x_e}\varphi_i - \partial_{x_e}\varphi_\ibar)_k + o(1/N)\Big)\Big\}\in[1,\infty),
\]
and observe that $a_N\to 1$ as $N\to\infty$. With this, using Jensen's inequality and the strong $\rmL^1$-convergence of $\iota_N(\mob_{i,\ibar,e}(c)) = K_{i\ibar} \sqrt{\rho_i \rho_\ibar \shift_e^N(\rho_i \rho_\ibar)}$ to $K_{i\ibar}\rho_i\rho_\ibar$, we obtain
\begin{align*}
    &\calR^\ast_{N,\swap} (c^N,\nablaDDiff\iota^*_N\varphi) 
    \le a_N \frac{1}{N^d} \sum_{i\in I} \sum_{\ibar\in\Ifull} \sum_{k\in\bbZ_N^d} \sum_{e\in \Eu} \mob_{i,\ibar,k,e}(c) \iota^*_N(|\partial_{x_e}\varphi_i - \partial_{x_e}\varphi_\ibar|_k^2) + o(1)\\
    &\qquad\qquad= a_N\sum_{i\in I} \sum_{\ibar\in\Ifull} \sum_{e\in \Eu}\big\langle \iota_N(\mob_{i,\ibar,e}(c)),\abs{\partial_{x_e}\varphi_i - \partial_{x_e}\varphi_\ibar}^2\big\rangle + o(1)
    \overset{N\to \infty}{\longrightarrow} \cR_\swap^*(\rho,\nablaCDiff\varphi).
\end{align*}
This concludes the proof.
\end{proof}

We are now in a position to deduce the liminf inequality for both $\calR_\swap$ and $\calR_\react$.
\begin{proposition}[Liminf for rates]\label{prop:liminf_R}
   Let $(c^N,F^N,J^N)$ satisfy the a priori estimates \eqref{eq:AprioriEstimBasic}. 
   Let $\rho$ be the limit of $(\iota_N c^N)_N$ from Propositions~\ref{prop:FluxBound}, \ref{prop:space_reg_bar-rho}, and \ref{prop:strong_compactness}. 
   Let $j$ be the limit of $\iota_{N,\react} J^N$ from Proposition \ref{prop:FluxBound}. Then, we have
    \[
    \liminf_{N\to\infty} \int_0^T \calR_N(c^N,F^N,J^N)\dd t\geq\int_0^T \mathcal{R}(\rho,f,j) \dd t.
    \]
\end{proposition}
\begin{proof}
    In light of Lemma~\ref{lem:EstimateDualDissipationPotential} and the assumption~\eqref{as_general:react} ensuring either the strong $\rmL^1$-convergence $(\rho^N)^{(\alpha^r+\beta^r)/2} \to \rho^{(\alpha^r+\beta^r)/2}$ if $r\in \Rone$, or the concavity of the map $\rho\mapsto\rho^{(\alpha^r+\beta^r)/2}$ if $r\in \Rtwo$, this follows exactly as in the proof of \cite[Proposition~5.14]{HMS25}.
\end{proof}

Finally, we are in a position to prove Theorem \ref{thm:Comp_Liminf}.
\begin{proof}[Proof of Theorem \ref{thm:Comp_Liminf}]
The asserted liminf inequality for the dissipations follows from Propositions~\ref{prop:liminf_S} and~\ref{prop:liminf_R} together with the compactness results Propositions \ref{prop:FluxBound} -- \ref{prop:space_reg_rho}.

Regarding the liminf inequality for the energies, we note that by definition of $\iota_N$ we have $\calE_N(c^N) = \calE(\rho^N)$. Therefore, as $\calE$ is an integral functional with convex integrand, the liminf inequality follows from the pointwise-in-time weak-$^\ast$ convergence of $\rho^N\rightharpoonup^\ast\rho$ established in Proposition \ref{prop:FluxBound}.
\end{proof}

\section{Limit chain rule}\label{s:CR}

The goal of this section is the proof of the chain rule identity Theorem~\ref{thm:CRContinuous} under the assumption that all species are fully connected.

\begin{proof}[Proof of Theorem~\ref{thm:CRContinuous}.]
Let us recall that by~\eqref{as_general:pos_coeff} we have $\Ifull = I$, so that by the bound $\calD(\rho,f,j)<\infty$ it holds $\nabla\sqrt{\rho_i}\in\rmL^2([0,T]\times \bbT^d)$ for all species $i\in I$. 
Our goal is to show for all $0\le s\le t\le T$ the identity
\begin{align}\label{eq:CR}
    \calE(\rho(t)) - \calE(\rho(s))
    &= \int_s^t\int_{\bbT^d}\bra[\Bigg]{\sum_{i,\ibar\in I} \bra[\big]{\nabla\rmD_i\calE(\rho) - \nabla\rmD_\ibar\calE(\rho)} \cdot f_{i,\ibar} 
    + \sum_{r\in R}(\Gamma\rmD\calE(\rho))_r j_r} \dd x \dd \tau.
\end{align}
To achieve this, we follow the classical strategy of regularizing and passing to the limit.
\paragraph{Step 1: Regularized chain rule.} We introduce the truncated energy
\begin{align*}
    \calE_\sigma(\rho) \coloneqq \sum_{i\in I}\int_{\bbT^d} \lambda_B(\rho_i + \sigma)\dd x
    \qquad\text{and hence}\qquad
    \rmD_i\calE_\sigma(\rho) = \log(\rho_i + \sigma).
\end{align*}
Next, we extend $(\rho,f,j)$ to $[-T,2T]$ by setting $\rho(t) = \rho(-t)$ for $t\in[-T,0]$, $\rho(t) = \rho(2T-t)$ for $t\in[T,2T]$, $(f,j)(t) = (-f,-j)(-t)$ on $[-T,0]$, and $(f,j)(t) = (-f,-j)(2T-t)$ on $[T,2T]$.
Then, we mollify them in both space and time and shift $\rho$ by $\overline\sigma =(\sigma,\ldots,\sigma)^\top \in\bbR^I$ for $\sigma\in(0,1]$, to obtain on $[0,T]$ the smooth triple $(\rho^\eps+\overline\sigma,f^\eps,j^\eps)$. 
By the linearity of the continuity equation, this triple is in $\CE^{\maxdensity+\abs{I}\sigma}$, so that using $(\log(\rho_i^\eps + \sigma))_i$ as test function in~\eqref{eq:CE_distributional} we obtain
\begin{align}\label{eq:CR-sigma-eps}
    \calE_\sigma(\rho^\eps(t)) {-} \calE_\sigma(\rho^\eps(s))
    &= \int_s^t\int_{\bbT^d}\bra[\Bigg]{\sum_{i,\ibar\in I} \bra[\big]{\nabla\rmD_i\calE_\sigma(\rho^\eps) {-} \nabla\rmD_\ibar\calE_\sigma(\rho^\eps)} \cdot f_{i,\ibar}^\eps 
    + \sum_{r\in R}(\Gamma\rmD\calE_\sigma(\rho^\eps))_r j_r^\eps} \dd x \dd \tau.
\end{align}
\paragraph{Step 2: $\eps\to 0$.}
Since, $\rho$ is narrowly continuous and $\calE_\sigma$ is narrowly lsc.\ and convex, so that $\calE_\sigma(\rho^\eps) \le \calE_\sigma(\rho)$, convergence of the left-hand side of \eqref{eq:CR-sigma-eps} follows for all $0\le s\le t \le T$ directly from dominated convergence and the bound $\rho\in \rmL^\infty([0,T]\times X)$. 
Next, we consider the first term on the right-hand side.
For all $i,\ibar\in I$ we have $f_{i,\ibar}^\eps \to f_{i,\ibar}$ strongly in $\rmL^2([0,T]\times\bbT^d)$ (e.g. bounding by the Hardy–Littlewood maximal function).
For the other factor we exploit the truncation to obtain for all $i\in I$ a.e.\ on $[0,T]\times\bbT^d$ that
\begin{align*}
    \abs{\nabla\rmD_i\calE_\sigma(\rho^\eps)} &=
    \abs{\nabla\log(\rho_i^\eps+\sigma)}
    \le \frac{1}{\sigma}\abs{\nabla\rho_i^\eps} \le \frac{1}{\sigma}\mathrm{M}(\nabla\rho_i),
\end{align*}
where we denoted by $\rmM(f)$ the Hardy–Littlewood maximal function of $f$. 
Recalling that $\abs{\nabla\rho_i} = 2\sqrt{\rho_i}\abs{\nabla\sqrt{\rho_i}}\in\rmL^2([0,T]\times \bbT^d)$, this yields for all $i\in I$ that $\nabla\rmD_i\calE_\sigma(\rho^\eps) \to \nabla\rmD_i\calE_\sigma(\rho)$ strongly in $\rmL^2([0,T]\times \bbT^d;\bbR^d)$ (and hence also the difference of two such terms) by dominated convergence.
Next, we argue similarly for the second term on the right-hand side. Indeed, recalling that $\rho_i\in[0,\maxdensity]$ for all $i\in I$ a.e.\ on $[0,T]\times\bbT^d$, we have for all $r\in R$ a.e.\ on $[0,T]\times\bbT^d$ the bound
\begin{align*}
    \abs{(\Gamma\rmD\calE_\sigma(\rho^\eps))_r j_r^\eps} 
    &\le \sum_{i\in I}\abs[\big]{\gamma^r_i\log(\rho_i^\eps + \sigma) j_r^\eps}
    \le \norm{\gamma^r}_1\bra[\big]{\log(\maxdensity+\sigma) - \log\sigma} \rmM(j_r),
\end{align*}
and the Hardy–Littlewood maximal function of $j_r$ is integrable since $\int_0^T\calR_\react(\rho,j)\dd t<\infty$, $\rho\in \rmL^\infty([0,T]\times X)$, and $\C$ is superlinear.
Therefore, $(\Gamma\rmD\calE_\sigma(\rho^\eps))_r j_r^\eps \to (\Gamma\rmD\calE_\sigma(\rho))_r j_r$ strongly in $\rmL^1([0,T]\times \bbT^d;\bbR^d)$ by dominated convergence, so that it holds
\begin{align}\label{eq:CR-sigma}
    \calE_\sigma(\rho(t)) {-} \calE_\sigma(\rho(s))
    &= \int_s^t\int_{\bbT^d}\bra[\Bigg]{\sum_{i,\ibar\in I} \bra[\big]{\nabla\rmD_i\calE_\sigma(\rho) {-} \nabla\rmD_\ibar\calE_\sigma(\rho)} \cdot f_{i,\ibar} 
    + \int_{\bbT^d}\sum_{r\in R}(\Gamma\rmD\calE_\sigma(\rho))_r j_r} \dd x\dd \tau.
\end{align}
\paragraph{Step 3: $\sigma\to 0$.}
Passing to the limit on the left-hand side of \eqref{eq:CR-sigma} is immediate for all times.
For the first term on the right-hand side the limit once more follows from dominated convergence, upon observing that for all $i,\ibar\in I$ and a.e.\ on $[0,T]\times\bbT^d$ we have the bound
\begin{align*}
    \abs{\nabla\rmD_i\calE_\sigma(\rho) \cdot f_{i,\ibar}} 
    = \abs{\nabla\log(\rho_i+\sigma)\cdot f_{i,\ibar}}
    = 2\abs[\bigg]{\nabla\sqrt{\rho_i+\sigma}\cdot \frac{f_{i,\ibar}}{\sqrt{\rho_i+\sigma}}}
    \le 2\abs{\nabla\sqrt{\rho_i}} \frac{\abs{f_{i,\ibar}}}{\sqrt{\rho_i}},
\end{align*}
and that the right-hand side term is in $\rmL^1([0,T]\times\bbT^d)$ due to boundedness of the dissipation.
For the second term we use Fenchel's inequality to obtain for all $\sigma\in(0,1]$, $r\in R$ a.e.\ on $[0,T]\times\bbT^d$ the bound
\begin{align*}
    \abs{(\Gamma\rmD\calE_\sigma(\rho))_r j_r} 
    &= \abs[\big]{\log\bra[\big]{(\rho + \overline\sigma)^{\gamma^r}}} \abs{j_r}\\
    &\le \rho^{(\alpha^r+\beta^r)/2}\rmC^\ast\bra[\big]{\log\bra[\big]{(\rho + \overline\sigma)^{\gamma^r}}}
    + \rmC\bra[\big]{j_r\big|\rho^{(\alpha^r+\beta^r)/2}}\\
    &= 2\rho^{(\alpha^r+\beta^r)/2}\bra[\Big]{(\rho + 1)^{(\alpha^r-\beta^r)/2} + (\rho + 1)^{(\beta^r-\alpha^r)/2}-2}
    + \rmC\bra[\big]{j_r\big|\rho^{(\alpha^r+\beta^r)/2}}\\
    &\le 2\bra[\Big]{(\rho + 1)^{\alpha^r} + (\rho + 1)^{\beta^r} - 2\rho^{(\alpha^r+\beta^r)/2}}
    + \rmC\bra[\big]{j_r\big|\rho^{(\alpha^r+\beta^r)/2}},
\end{align*}
where we denoted by $1$ also the vector in $\bbR^I$.
The dissipation bound also ensuring integrability of this right-hand side, one more application of the dominated convergence theorem yields~\eqref{eq:CR} for a.e.\ $0\le s\le t\le T$. 

\paragraph{Step 4: Equivalence of solution concepts}

Combining~\eqref{eq:CR} with Fenchel's inequality, we obtain for all $0\le s \le t \le T$
\begin{align*}
    \calE(\rho(t)) - \calE(\rho(s)) &= \int_s^t\int_{\bbT^d}\bra[\Bigg]{\sum_{i,\ibar\in I} 2\bra[\big]{\sqrt{\rho_\ibar}\nabla\sqrt{\rho_i} - \sqrt{\rho_i}\nabla\sqrt{\rho_\ibar}} \cdot \frac{f_{i,\ibar}}{\sqrt{\rho_i\rho_\ibar}}
    + \sum_{r\in R}\log\bra[\big]{\rho^{\gamma^r}} j_r} \dd x \dd \tau\\
    &\ge -\int_s^t\int_{\bbT^d}\sum_{i,\ibar\in I} \bra[\bigg]{2K_{i\ibar}\abs[\big]{\sqrt{\rho_\ibar}\nabla\sqrt{\rho_i} - \sqrt{\rho_i}\nabla\sqrt{\rho_\ibar}}^2 + \frac{\abs{f_{i,\ibar}}^2}{2K_{i\ibar}\rho_i\rho_\ibar} }\dd x \dd \tau\\
    &\quad-\int_s^t\int_{\bbT^d}\sum_{r\in R}\bra[\Big]{2\kappa_r \abs[\big]{\rho^{\alpha^r/2} - \rho^{\beta^r/2}}^2 + \rmC\bra[\big]{j_r\big|\kappa_r\rho^{(\alpha^r+\beta^r)/2}}}\dd x \dd \tau\\
    &= -\calD^{[s,t]}(\rho,f,j),
\end{align*}
with equality if and only if it holds the constitutive relation~\eqref{eq:const_rel_cont}.
In turn,~\eqref{eq:const_rel_cont} and the continuity equation~\eqref{eq:CE_distributional} are together equivalent to the PDE~\eqref{eq:reaction-cross-diffusion_intro}, whence the claimed equivalence of the two solution concepts.
\end{proof}

\begin{remark}
    The above proof exploits the bound $\sqrt{\rho_i} \in \rmL^2([0,T];\rmH^1(\bbT^d))\cap \rmL^\infty([0,T]\times\bbT^d)$ to control both the shift and the mollification. 
    Without~\eqref{as_general:pos_coeff}, we only have $\rmL^2$-control over terms of the form $\sqrt{\rho_\ibar}\,\nabla\sqrt{\rho_i}$, i.e., 
    we lose control over $\nabla\sqrt{\rho_i}$ where $\rho_\ibar$ vanishes. 
    Mollifying the two factors individually, decouples them and thus destroys the upper bounds that would allow for an application of dominated convergence or similar arguments.
    
    Regarding the reactions, note that in \cite{HMS25} it was shown, that $\rmL^\infty$-bounds are not needed to prove of the chain rule for a class of stoichiometric coefficients beyond what implies convexity, but more restrictive than what we consider here.
\end{remark}

\section{Tilting the energy landscape}\label{s:tilting}

In this section, we discuss the adaptation of the previous steps to include drift potentials $V\in \rmL^\infty(X)\cap\rmH^1(X)$ with $\sum_{i\in I}\exp(-V_i) = M$ a.e.\ on $\bbT^d$ in the analysis, thereby proving Theorem~\ref{thm:weighted_main-results}.
To this end, we introduce a reference concentration $\omega\in \rmL^\infty(X)\cap\rmH^1(X)$ given by the identities $\omega_i \coloneqq \exp(-V_i)$ for all $i\in I$ and hence satisfies $\omega_\ast \le \omega_i\le \omega^\ast$ for all $i\in I$ and some positive finite constants $\omega_\ast,\omega^\ast\in (0,\infty)$.

\paragraph{Compactness.}
To translate the results from Section~\ref{s:Compactness} we need to ensure that the uniform bounds~\eqref{eq:AprioriRegularity} remain true in the weighted setting. 
Then, all arguments from that section remain valid.
To show the uniform bounds, for $N\in\bbN$, $a:\bbZ^d_N\to \R$, $k\in\bbZ^d_N$ and $e\in\Eu$, recall the discrete gradient $\babla a_{k,e} = a_{k+e}-a_k$ and the arithmetic mean $\theta(a)_{k,e} = (a_k+a_{k+e})/2$, and denote the geometric mean by $\tilde\theta(a)_{k,e} =  \sqrt{a_k a_{k+e}}$.
With these notations we state the following lemma.
\begin{lemma}\label{lem:discrete_calculations_pointwise}
    Let $c,C>0$ and consider $a, \bar a:\bbZ^d_N\to  [0,C]$ and $b,\bar b:\bbZ^d_N\to [c,C]$.
    Then, on $\bbZ^d_N\times\Eu$ we have the inequality 
    \begin{equation}
    \label{eq:est-relativ_slope-absolute_slope}
    \begin{aligned}
        \tilde\theta(b^2)\tilde\theta(\bar b^2)\abs[\bigg]{\theta\bra[\bigg]{\frac{a}{b}}\babla\bra[\bigg]{\frac{\bar a}{\bar b}} - \theta\bra[\bigg]{\frac{\bar a}{\bar b}}\babla\bra[\bigg]{\frac{a}{b}}}^2
        \ge \frac12\abs[\bigg]{\theta(a)\babla \bar a - \theta(\bar a)\babla a }^2 -  \bar C(\abs{\babla b}^2 + \abs{\babla \bar b}^2),
    \end{aligned}
    \end{equation}
    where $\bar C>0$ only depends on $c,C>0$.
\end{lemma}
\begin{proof}
    We first make some elementary observations.
    Given $a:\bbZ^d_N\to  [0,C]$ and $b:\bbZ^d_N\to [c,C]$, we have on $\bbZ^d_N\times\Eu$ the chain-rule type identity $\babla b^{-1} = -\tilde\theta\bra[\big]{b^{-2}}\babla b$.
    With this, a short calculation and the estimate $|\babla \sqrt{b}|^2 \le |\babla b|$ yield
    \begin{align}
        \notag
        \tilde\theta(b) \theta\bra[\big]{b^{-1}} 
        &= 1 + \frac12 \tilde\theta\bra[\big]{b^{-1}}\abs[\big]{\babla\sqrt{b}}^2
        = 1+ O(\babla b),\\
        \tilde\theta(b) \theta\bra[\big]{b^{-1}a}
        \notag
        &= \tilde\theta(b) \theta\bra[\big]{b^{-1}}\theta(a) + \frac14 \tilde\theta(b)\babla b^{-1}\babla a
        = \theta(a)\pra[\Big]{1 + \frac12 \tilde\theta\bra[\big]{b^{-1}}\abs[\big]{\babla\sqrt{b}}^2} -\frac14 \tilde\theta\bra[\big]{b^{-1}} \babla b\babla a\\
        \notag
        &= 
        \theta(a) + O(\babla b).
    \end{align}
    Using these estimates, the discrete product identity $ \babla \bra[\big]{ab} = \theta(a)\babla b + \theta(b)\babla a$, and Young's inequality, we then obtain \eqref{eq:est-relativ_slope-absolute_slope} as follows.

\begin{align*}
    &\phantom{==}\tilde\theta(b^2)\tilde\theta(\bar b^2)\abs[\bigg]{\theta\bra[\bigg]{\frac{a}{b}}\babla\bra[\bigg]{\frac{\bar a}{\bar b}} - \theta\bra[\bigg]{\frac{\bar a}{\bar b}}\babla\bra[\bigg]{\frac{a}{b}}}^2\\
    &=
    \tilde\theta(b^2)\tilde\theta(\bar b^2)\abs[\bigg]{\theta\bra[\bigg]{\frac{a}{b}}\pra[\Big]{\theta(\bar b^{-1})\babla \bar a + \theta(\bar a)\babla \bar b^{-1}} - \theta\bra[\bigg]{\frac{\bar a}{\bar b}}\pra[\Big]{\theta(b^{-1})\babla a + \theta(a)\babla b^{-1}} }^2\\
    &=
    \abs[\bigg]{\tilde\theta(b)\theta\bra[\bigg]{\frac{a}{b}}\pra[\Big]{\tilde\theta(\bar b)\theta(\bar b^{-1})\babla \bar a + \tilde\theta(\bar b)\babla \bar b^{-1}\theta(\bar a)} - \tilde\theta(\bar b)\theta\bra[\bigg]{\frac{\bar a}{\bar b}}\pra[\Big]{\tilde\theta(b)\theta(b^{-1})\babla a + \tilde\theta(b)\babla b^{-1}}\theta(a) }^2\\
    &\ge \frac12\abs[\big]{\theta(a)\babla \bar a - \theta(\bar a)\babla a }^2 -  \bar C\bra[\big]{\abs[\big]{\babla b}^2 + \abs[\big]{\babla \bar b}^2}.
\end{align*}
This concludes the proof.
\end{proof}
With the above lemma, we are now in a position to establish the required bound.
\begin{proposition}
    Let~\eqref{eq:AprioriEstimBasic} hold true for the modified energy and slope functionals introduced in Theorem~\ref{thm:weighted_main-results}.
    Then,~\eqref{eq:AprioriRegularity} holds as well.
\end{proposition}
\begin{proof}
    Fixing and dropping an index $i\in I$ and arguing as in the proof of Lemma~\ref{lem:CE.DiscCont}, we have
    \begin{align*}
        \frac{1}{N^d}\sum_{k\in\bbZ^d_N}\sum_{e\in\Eu} N^2\abs[\big]{\babla\sqrt{w^N}_{k,e}}^2 
        &= \frac{1}{N^d}\sum_{k\in\bbZ^d_N}\sum_{e\in\Eu} \frac{N^2 \abs{\babla w^N_{k,e}}^2}{4(\theta(\sqrt{w^N})_{k,e})^2}\\
        &\le \frac{1}{4\omega_\ast} \sum_{e\in\Eu}\int_{\bbT^d}\abs[\bigg]{\int_0^1 \shift^N_{e\theta}\partial_{x_e}\omega(x) \dd \theta}^2\dd x
        \le \frac{1}{4\omega_\ast} \int_{\bbT^d}\abs[\big]{\nabla\omega(x)}^2\dd x,
    \end{align*}
    where for the last inequality we used H{\"o}lder's inequality and the Fubini-Tonelli theorem.
    Therefore, the error terms arising on the right-hand side of~\eqref{eq:est-relativ_slope-absolute_slope} when applying Lemma~\ref{lem:discrete_calculations_pointwise} for $a =\sqrt{c}$ and $b=\sqrt{w^N}$ are controlled so that the result follows once more from Lemma~\ref{lem:discrete_calculations_sum}.
\end{proof}

\paragraph{Lower limit inequalities.}
The considerations about the rate term $\calR$ in Section~\ref{s:Liminf} remain unchanged and
the convergence of the energies follows from the strong $\rmL^p$-convergence $\omega^N\to\omega$ for all $p\in[1,\infty)$ and the arguments for the non-tilted case.

It remains to consider the slope $\calS$, which now contains quotients with the reference measure.
This on its own is not a complication.
However, to obtain the lower limit in $\calS_\swap$ we now also have to accommodate for the different interpolations used on $\omega$.
Here, simply combining the established limits and controlling arising commutator errors is not optimal as controlling the commutator errors in $\rmL^\infty$ would require continuity of $\omega$.
We avoid this by directly proving the weak $\rmL^2$-convergence of the integrand including the different mean functions.
\begin{proposition}[Improved liminf for slopes]\label{prop:liminf_S_omega}
    Let $\omega$, $\calS$, $\discref^N$, $\calS_N$ be as in Theorem~\ref{thm:weighted_main-results}.
    Let $\rho^N = \iota_Nc^N \to \rho$ strongly in $\rmL^2([0,T]\times X)$ and satisfy the gradient limit~\eqref{eq:babla_conv}.
    Then, it holds
    \begin{align*}
        \liminf_{N\to\infty} \calS_N(c^N) \ge \calS(\rho).
    \end{align*}
\end{proposition}
\begin{proof}
    We first consider $\calS_\swap$.
    Denoting $\omega^N = \iota_N\discref^N = \iota_N\iota_N^\ast \omega$, for all $(i,\ibar)\in I{\times}\Ifull\cup \Ifull{\times}I$, and $e\in\Eu$, (along a subsequence) it holds weakly in $\rmL^2([0,T]\times \bbT^d)$ the limit
    \begin{align*}
        \frac N2\sqrt{\shift^N_e\bra[\big]{\sqrt{\omega^N_\idot}} \sqrt{\omega^N_\idot} \shift^N_e\bra[\big]{\sqrt{\omega^N_\ibar}} \sqrt{\omega^N_\ibar}} \bra[\Big]{\shift^N_e\bra[\Big]{\sqrt{\tfrac{\rho^N_\idot}{\omega^N_\idot}}}+ \sqrt{\tfrac{\rho^N_\idot}{\omega^N_\idot}}} \babla_N\sqrt{\smash{\tfrac{\rho^N_\ibar}{\omega^N_\ibar}}\vphantom{\tfrac{\rho_\ibar}{\omega_\ibar}}}\one_{\set{\rho_\ibar>0}} \rightharpoonup \sqrt{\rho_i\omega_\ibar}\nabla\sqrt{\tfrac{\rho_\ibar}{\omega_\ibar}}.
    \end{align*}
    Indeed, this is a straight-forward adaptation of Proposition~\ref{prop:space_reg_rho}, where existence of a weak $\rmL^2$-limit of the left-hand side follows from~\eqref{eq:AprioriRegularity} and the regularity of $\omega$, and identifying its weak $\rmL^1$-limit (which coincides with the weak $\rmL^2$-limit) as the right-hand side is possible using the strong $\rmL^2$-convergence $\omega^N\to \omega$.
    With the limit established, the liminf inequality follows as in the proof of Proposition~\ref{prop:liminf_S}.

    Finally, note that for $\calS_\react$ there is no spatial interpolation involved and hence, due to the bounds on $\omega$, no additional complications arise when including it.
    This completes the proof.
\end{proof}

\paragraph{Chain rule identity.}
Due to the bounds and $\rmH^1$-regularity of $\omega$, it does not introduce new complications so that the presented arguments carry over to the more general setting. 

Alternatively, it is straight forward to first prove the chain rule identity for the potential energy, alone, and then add the previously established untilted chain rule identity.

\bibliographystyle{alphaabbr}
{\small \bibliography{bibliography}}

\end{document}